\documentclass[11pt]{amsart}

\usepackage{amsmath, amsfonts, amssymb, amscd, verbatim, xypic, xy}
\usepackage{pstricks, pst-node, pst-tree}

\newcommand{\marginlabel}[1]%
  {\mbox{}\marginpar{\raggedleft\hspace{0pt}\bfseries\sf#1}}

\def\ZZ{{\mathbb Z}}
\def\NN{{\mathbb N}}

\def\CC{{\mathbb C}}

\def\QQ{{\mathbb Q}}
\def\PP{{\mathbb P}}

\def\cC{\mathcal{C}}

\def\cA{\mathcal{A}}

\def\cG{\mathcal{G}}

\def\cL{\mathcal{L}}
\def\cO{\mathcal{O}}
\def\cM{\mathcal{M}}
\def\cN{\mathcal{N}}
\def\cP{\mathcal{P}}
\def\cR{\mathcal{R}}

 \DeclareMathOperator{\Coker}{coker}

\DeclareMathOperator{\sym}{sym}

\newtheorem{lemma}{Lemma}[section]
\newtheorem{theorem}[lemma]{Theorem}
\newtheorem{corollary}[lemma]{Corollary}
\newtheorem{proposition}[lemma]{Proposition}
\newtheorem*{claim}{Claim}

\theoremstyle{definition}
\newtheorem{definition}[lemma]{Definition}
\newtheorem{example}[lemma]{Example}

\newtheorem{remark}[lemma]{Remark}

\numberwithin{equation}{section}

\newcommand{\bean}{\begin{eqnarray}}
\newcommand{\eean}{\end{eqnarray}}
\newcommand{\be}{\begin{displaymath}}
\newcommand{\ee}{\end{displaymath}}
\newcommand{\bea}{\begin{eqnarray*}}
\newcommand{\eea}{\end{eqnarray*}}

\newcommand{\ol}{\overline}

\psset{levelsep=20pt,nodesep=0pt,treesep=10pt}

\begin{document}

\title{The Chow ring of $\ol{M}_{0,m}(\PP^n, d)$}

\author[Anca M. Musta\c{t}\v{a}]{Anca~M.~Musta\c{t}\v{a}}
\author[Andrei Musta\c{t}\v{a}]{Andrei~Musta\c{t}\v{a}}
\address{Department of Mathematics, University of Illinois,
Urbana-Champaign, IL, USA} \email{{\tt amustata@math.uiuc.edu,
dmustata@math.uiuc.edu}}

\date{\today}

\begin{abstract}
 We describe the Chow ring
    with rational coefficients of the moduli space of stable maps with marked points $\ol{M}_{0,m}(\PP^n,d)$ as
the subring of invariants of a ring
    $B^*(\ol{M}_{0,m}(\PP^n,d);\QQ) $, relative to the action of the
    group of symmetries $S_d$. 
   $B^*(\ol{M}_{0,m}(\PP^n,d);\QQ)$ is computed by following
 a sequence of intermediate spaces for $\ol{M}_{0,m}(\PP^n,d)$ and relating them to substrata of $\ol{M}_{0,1}(\PP^n,d+m-1)$. An additive basis for $A^*(\ol{M}_{0,m}(\PP^n,d);\QQ) $ is given. 
\end{abstract}
\maketitle

\bigskip

\section*{Introduction}

The moduli spaces of stable maps introduced by Kontsevich and Manin in \cite{kontsevich} have come to play a central role in the enumerative geometry of curves. In this paper we are concerned with $\ol{M}_{0,m}(\PP^n,d)$ when $m>0$, which parametrizes stable maps from genus zero, at most  nodal curves with $m$ distinct, smooth marked points, into the projective space. $\ol{M}_{0,m}(\PP^n,d)$ is a smooth Deligne-Mumford stack having a projective coarse model (see \cite{behrend1}, \cite{fultonpan}). Thus the fine and coarse moduli spaces share the same cohomology and Chow groups  with rational coefficients (\cite{vistoli}). Intersection numbers on  $\ol{M}_{0,m}(\PP^n,d)$, known as genus zero Gromov-Witten invariants, have been intensively studied starting with \cite{kontsevich}, \cite{givental}, \cite{pandharipande2}. 

 Recent results regard other aspects of the cohomology: the cohomology groups are generated by tautological classes (\cite{oprea1}). An algorithm for computing the Betti numbers of all spaces $\ol{M}_{0,m}(\PP^n,d)$ is given in \cite{getzlerpan}. Forays into the cohomology ring structure for specific cases of moduli spaces of stable maps start with  the degree zero case $\ol{M}_{0,m}$, whose Chow ring is calculated in \cite{keel}.
  The cohomology of the moduli space without marked points $\ol{M}_{0,0}(\PP^n,d)$ was computed in degree 2 and partially in degree 3 by \cite{behrend2}. An additive basis for $\ol{M}_{0,2}(\PP^n,2)$ was written down in \cite{cox}, the ring structure of $\ol{M}_{0,2}(\PP^1,2)$ was determined in \cite{cox2}. \cite{noi} contains a description of  the Chow ring
 with rational coefficients of the moduli space  $\ol{M}_{0,1}(\PP^n,d)$ as the subalgebra of invariants of a larger $\QQ$-algebra
$B^*(\ol{M}_{0,1}(\PP^n,d);\QQ) $, relative to the action of the group of symmetries $S_d$. The study of  the algebra $B^*(\ol{M}_{0,1}(\PP^n,d);\QQ) $ is both motivated and made possible by the existence of a sequence of intermediate spaces for $\ol{M}_{0,1}(\PP^n,d)$, the birational morphisms between which can be understood in detail.

 In the present paper we extend this construction to $ \ol{M}_{g,m}(\PP^n,d)$, whose intermediate moduli spaces turn out to be intimately related to substrata of $\ol{M}_{g,1}(\PP^n,d+m-1)$. In genus 0, this leads to a natural construction of the ring $B^*(\ol{M}_{0,m}(\PP^n,d);\QQ)$, whose structure, presented in Theorem 5.1, follows essentially from work done in \cite{noi}.  An additive basis for $A^*(\ol{M}_{0,m}(\PP^n,d);\QQ) $ is also constructed (Proposition 6.1). The elements of the basis are clearly tautological and can be counted as a sum indexed by rooted decorated trees with marked leaves. 

  The paper is organized as follows: The first section is devoted to the construction of a system of intermediate spaces for $\ol{M}_{g,m}(\PP^n,d)$. They are Deligne-Mumford stacks finely representing certain moduli problems, and in genus 0 they are smooth. Via the rigidification introduced in \cite{fultonpan}, these moduli spaces turn out to be closely related to moduli spaces of rational weighted stable curves previously studied by \cite{hassett}. Indeed, rigidified moduli spaces form a finite  cover of $\ol{M}_{g,m}(\PP^n,d)$ and  are at the same time $(\CC^*)^n$-- torus bundles over locally closed subsets of $ \ol{M}_{g,m+(n+1)d}$. This is visible when taking into account, for each generic point $(C,\{p_i\}_{1\leq i\leq m}\to \PP^n)$ of $\ol{M}_{g,m}(\PP^n,d)$,  $(n+1)$ independent hyperplane sections on $C$. On the other hand, a system of weights on the marked points allows Hassett in \cite{hassett} to construct moduli spaces of curves where some marked points are allowed to coincide.  Quotients by $(S_d)^n$ of torus bundles over subspaces of Hassett's moduli spaces glue together to stacks birational to $\ol{M}_{g,m}(\PP^n,d)$. These are the intermediate (weighted) moduli spaces of stable maps  $\ol{M}_{g,\cA}(\PP^n,d,a)$. They parametrize rational curves with weighted marked points mapping to $\PP^n$, such that among their rational tails, only those of certain degrees are allowed.  

 The second and third sections are dedicated to keeping record of the canonical stratification of $\ol{M}_{0,m}(\PP^n,d)$ and its intermediate spaces by means of decorated trees and good monomials of stable two-partitions. The sources of inspiration for this approach are  \cite{behrend1}, \cite{km}, \cite{losev}. Some extra structure is attached to the trees to account both for the substrata of $\ol{M}_{0,\cA}(\PP^n,d,a)$, and for the complex network of morphisms between them. In this context the language of 2-partitions fits well with the construction of the $\QQ$-algebra   $B^*(\ol{M}_{0,\cA}(\PP^n,d,a))$.

  The ring  $B^*(\ol{M}_{0,\cA}(\PP^n,d,a))$ introduced in section 4 holds inside it the Chow groups of all substrata of $ \ol{M}_{0,\cA}(\PP^n,d,a)$, and takes into account all morphisms between substrata indexed by contractions of trees. There is a natural action  of the group of permutations $S_d$ on the generators of $B^*(\ol{M}_{0,\cA}(\PP^n,d,a))$, and the invariant subalgebra with respect to this action is  $A^*(\ol{M}_{0,\cA}(\PP^n,d,a))$.

 A formula for   $B^*(\ol{M}_{0, m}(\PP^n,d))$ is given in Theorem 5.1, and an additive basis for $A^*(\ol{M}_{0, m}(\PP^n,d))$ is written in Proposition 6.1. The ensuing Poincair\'e polynomial computation, as a sum indexed by rooted decorated trees with marked leaves, is especially practical in low-degree examples like those considered at the end of the paper.

We are grateful to Kai Behrend and Daniel Krashen for helpful suggestions, to Jim Bryan for interest and hospitality at the University of British Columbia. 
We thank Sheldon Katz and the referee for suggestions in improving the exposition. We thank Michael Guy for useful comments. 

 After writing this paper we were informed that Arend Bayer and Yuri Manin are also constructing moduli spaces of weighted pointed stable maps for smooth projective targets, based on Hassett's weighted curve spaces (\cite{bayer-manin}) and  motivated by interest in wall-crossing phenomena for different stability conditions.  These spaces have a good obstruction theory that leads to the construction of a virtual fundamental class.
We work in a complementary direction by considering weights on maps as well as on the marked points of a curve. This condition proves central in the cohomological study of stable map spaces (\cite{noi}, \cite{noi3}).  We would like to thank Prof. Yuri Manin for providing a copy of the preprint and explaining their motivation.

\section{Intermediate moduli spaces of $\ol{M}_{g,m}(\PP^n,d)$}

Let $g$ and $m$ be two integers, $g\geq 0$ and $m\geq 0$. 
In this section we consider complex curves of any arithmetic genus $g$,  with $m$ marked points and at most nodal singularities different from the markings, and morphisms to $\PP^n$.

Fix a rational number $a>0$ and an $m$-tuple $\cA=(a_1,...,a_m)\in \QQ^m$  of rational numbers  such that $0\leq a_j\leq 1$ for all $j=1,...,m$ and such that $\sum_{i=1}^{m}a_i+da>2-2g$. We will work over a field $k$ such that $d< \mbox{ char } k$ or $k$ has characteristic zero (as in \cite{behrend1}, Theorem 3.14).

\begin{definition}

 An $ (\cA, a)$-stable  family of degree $d$ nodal maps  from rational curves with $m$ marked sections to $\PP^n$  consists of  the following
 data:
\[( \pi \colon C \to S , \{p_i\}_{1\leq i\leq m}, \cL , e ) \]
where  $\cL$ is a line bundle on $C$ of degree $d$ on each fiber $C_s$, and $e:\cO^{n+1}\to \cL$ is a morphism of fiber bundles (specified up to isomorphisms of the target) such that:

\begin{enumerate}
\item $\omega_{C|S}(\sum_{i=1}^ma_ip_i)\otimes \cL^a$ is relatively ample over $S$,
\item $\cG :=\Coker e$, restricted over each fiber $C_s$, is a skyscraper sheaf supported only on smooth points, and 
\item for any $p \in S$ and for any $ I \subseteq \{1,...,m\}$ (possibly empty) such that $ p = p_i $ for all $i \in I$ the following condition holds  $$\sum_{\i \in I}  a_i + a\dim\cG_{p}\leq 1.$$ 

\end{enumerate}

\end{definition}

 The rigidified version of this moduli problem is based on Hassett's moduli spaces of curves. The following notions were introduced in \cite{hassett}.

\begin{definition}
 A collection of $N$ weights is an $N$-tuple $\cA=(a_1,...,a_N)\in \QQ^N$ such that $0<a_j\leq 1$ for all j=1,...,N and such that $\sum_{i=1}^{N}a_i>2-2g$, where $g \geq 0$.

 Given a collection of $N$ weights $\cA$, a family of genus $g$ nodal curves with marked points $\pi \colon (C, s_1,...,s_N)\to S$ is $(g,\cA)$-stable if:
\begin{enumerate}
\item The map $\pi: C\to S$ is flat, projective and such that each geometric fiber is a  nodal connected curve of arithmetic genus $g$.
\item The sections $(s_1,...,s_N)$ of $\pi$ lie in the smooth locus of $\pi$, and for any subset $\{s_i\}_{i\in I}$ with nonempty intersection, $\sum_{i\in I} a_i\leq 1$.
\item $K_{\pi} +\sum_{i=1}^N a_i s_i $ is relatively ample, where $K_{\pi}$ is  the relative canonical divisor.
\end{enumerate}
   \end{definition}

 \begin{theorem} (2.1 in \cite{hassett} ) The moduli problem of $(g,\cA)$-stable curves is finely represented by a smooth Deligne-Mumford stack $\ol{M}_{g,\cA}$, proper over $\ZZ$. The corresponding coarse moduli scheme is projective over $\ZZ$.

\end{theorem}

 A natural collection of morphisms between these spaces can be constructed:
\begin{theorem} (4.1 in \cite{hassett} )
For two collections of weights $\cA=(a_1,...,a_N)$ and
$(\cA'=(a'_1,...,a'_N))$ such that $a_i\leq a'_i$ for each $i$,
there exists a natural birational reduction morphism $\rho:
\ol{M}_{g,\cA'} \to  \ol{M}_{g,\cA}.$ The image of an element
$(C, s_1,...,s_N)\in\ol{M}_{g,\cA'} $ is obtained by successively
collapsing components along which the divisor $K_C+\sum_{i=1}^N
a_is_i$ fails to be ample.

\end{theorem}

Similarly, there are smooth Deligne-Mumford stacks of weighted stable curves with marked divisors, obtained as quotients of the above by groups of permutations of marked points.
Rigid  $ (\cA, a) $-stable maps are defined as follows:

\begin{definition}
Let  $\bar{t}=(t_0:...:t_n)$ denote a homogeneous coordinate
system on $\PP^n$. As before, let $a>0$ and  $\cA=(a_1,...,a_m)\in \QQ^m$ be an $m$-tuple of rational numbers  such that $0\leq a_j\leq 1$ for all $j=1,...,m$ and such that $\sum_{i=1}^{m}a_i+da>2-2g$.

A $ \bar{t}$-rigid, $ (\cA, a) $-stable family of degree $d$ nodal maps with $m$ marked points to $\PP^n$ consists of the following data:
\bea (\pi \colon C\to B,  \{q_{i,j} \}_{0 \leq i \leq n, 1 \leq j
\leq d}, \{p_i\}_{1\leq i\leq m}, \cL, e) \eea such that
\begin{enumerate}
\item  the family $(\pi \colon C \to S ,  \{
q_{i,j} \}_{0 \leq i \leq n, 1 \leq j \leq d}, \{p_i\}_{1\leq
i\leq m} )$ is a $(g,\cA')$- stable family of curves, where the
system of weights $\cA'$ consists of: \begin{itemize}
\item $a_{i,j}=\frac{a}{(n+1)}$ for the sections $ \{q_{i,j} \}_{0 \leq i \leq n, 1 \leq j \leq d}$;
\item $a_i$ for $p_i$, $i=1,...,m$.
\end{itemize}
\item  $\cL$ is a line bundle on $C$ of degree $d$ on each fiber of $\pi$, and  $e:\cO_{C}^{n+1}\to\cL$ is
a morphism of sheaves  such
that, via the natural isomorphism $H^0(\PP^n, \cO_{\PP^n}(1))
\cong H^0(C,\cO_{C}^{n+1})$, there is an equality of Cartier
divisors
$$ ( e( \bar{t}_i)=0 ) = \sum_{j=1}^d q_{i,j}. $$
\end{enumerate}
\end{definition}

A comment on the possible choices of $a$ is in order. In Definition 1.2, the requirement that $a_{j}\leq 1$ in conjunction with condition (3) insures that a $(g, \cA)$--stable curve has a finite group of automorphisms fixing its marked points. In the case of $\bar{t}$-rigid maps, no upper bound need be placed on the weight of the marked points coming from the coordinate hyperplanes. Indeed, if any such point is on a component of a  $\bar{t}$-rigid curve, then other $n$ such points are automatically on the same component, which implies a finite automorphism group whenever $n>0$, $d>0$.   

We note that Definition 1.5. and the following proposition admit natural variations in which the sections $q_{i,j}$ are replaced by divisors $D_{i,k}$ such that $\sum_k  D_{i,k}= ( e( \bar{t}_i) )$. 

With the notations from Definition 1.5, the following proposition holds.

\begin{proposition}
The moduli problem of $ \bar{t}$-rigid, $ (\cA, a) $-stable curves
is represented by a quasiprojective Deligne-Mumford stack
$\ol{M}_{g,\cA}(\PP^n,d,a,\bar{t})$, the total space of a
$(\CC^*)^n$--bundle over a locally closed subset of $\ol{M}_{g,\cA'}$. In particular, 
$\ol{M}_{0,\cA}(\PP^n,d,a,\bar{t})$ is a quasiprojective scheme mapping onto an open subset of $\ol{M}_{0,\cA'}$.

\end{proposition}

The proof is identical to that of Proposition 3 in
\cite{fultonpan}. The total space of the torus bundle comes into
play because for any map $C\to \PP^n$, the pull-backs of the
hyperplane divisors $(\bar{t}_i=0)$ determine the map only up to
the action of $(\CC^*)^n$ on $ \PP^n$.

The group $(S_d)^{n+1}$ has a natural action on $\ol{M}_{g,\cA}(\PP^n,d,a,\bar{t})$. The resulting quotients for various coordinate systems $\bar{t}$ glue to a Deligne-Mumford stack: the stack of  $ (\cA, a) $-stable degree $d$ maps.

\begin{proposition}
The moduli problem of genus $g$, $ (\cA, a) $-stable degree $d$ nodal  maps
with $m$ marked points into $\PP^n$  is finely represented by a proper
 Deligne-Mumford stack $\ol{M}_{g,\cA}(\PP^n,d,a)$. When $g=0$, the stack is smooth.
\end{proposition}

\begin{proof}

 Here we introduce the short notation $\ol{M}$ for the category fibered by groupoids associated to  $\ol{M}_{g,\cA}(\PP^n,d,a)$,  and $\ol{M}(\bar{t})$ for  $\ol{M}_{g,\cA}(\PP^n,d,a,\bar{t})$.  $\ol{M}$ is constructed by gluing quotients of   $\ol{M}(\bar{t})$, as in \cite{fultonpan}, Proposition 4.

The fact that $\ol{M}$ is a stack follows from Grothendieck descent theory like with the usual moduli spaces of stable maps. We will prove the existence of an \'etale cover of $\ol{M}$ by a smooth Deligne-Mumford stack $U$. In the process we observe that the diagonal $\ol{M}\to \ol{M}\times\ol{M}$ is representable, finite type and separated. In genus 0, $U$ is smooth. The stack $\ol{M}$ is proper, as $\ol{M}_{g,n}(\PP^n,d)$ maps onto it.

Let $D=\{1,...,d\}$ and let $\cP_d$ denote the set of partitions of $D$. For each partition $b=\{B_1,...,B_l\}\in \cP_d$ we denote by $S_{b}$ the direct sum of symmetric groups $S_{|B_1|}\oplus ...\oplus S_{|B_l|}$, each $S_{|B_i|}$ acting on the subset $B_i$, and by $\mu_b:=\{|B_1|,...,|B_l|\}$ the associated partition of $d$. We will say that a partition $b'=\{B'_1,...,B'_k\}$ is a refinement of another partition $b=\{B_1,...,B_l\}$ if each set $B_i$ is a union of sets $B'_j$. Consider an element $\beta =(b^0,...,b^{n})$ in the $(n+1)$-th product $ \cP^{n+1}_d$ and let $\mu_{\beta }= (\mu_{b^0},...,\mu_{b^{n}})$. Let $l^i$ denote the length of the partition $b^i$. Thus $b^i=\{B^i_1,...,B^i_{l^i}\}$.

For $\mu_{\beta }$ and
 any homogeneous coordinate system $\bar{t}=(t_0:...:t_n)$ on $\PP^n$,
there is a moduli problem which we will denote by $\ol{\cM}(\mu_{\beta }, \bar{t})$, defined as follows. For any scheme $S$,  the set $\ol{\cM}(\mu_{\beta },\bar{t} )(S)$ consists of isomorphism classes of objects $( \pi \colon C \to S , \{D_{i,j}\}_{0\leq i\leq n, 1\leq j\leq l^i}, \{p_i\}_{1\leq i\leq m}, \cL , e )$ such that
\begin{enumerate}
\item for each $s\in S$, $D_{i,j}|_{C_s}$ is an effective Cartier divisor of degree $|B^i_j|$, disjoint from any other divisor $D_{i,j'}|_{C_s}$. The family $$(\pi \colon C \to S ,  \{D_{i,j}\}_{0\leq i\leq n, 1\leq j\leq l^i}, \{p_i\}_{1\leq
i\leq m} )$$ is an $(\cA'')$- stable family of curves with marked  points $\{p_i\}_{1\leq
i\leq m}$ and divisors $\{D_{i,j}\}_{0\leq i\leq n, 1\leq j\leq l^i}$, where  the
system of weights $\cA''$ consists of
 \begin{itemize} \item $a_{i,j}=\frac{a |B^i_j|}{(n+1)}$ for the divisors $\{D_{i,j}\}_{0\leq i\leq n, 1\leq j\leq l^i}$, \item $a_i$ for the sections $p_i$.\end{itemize}
\item  $\cL$ is a line bundle on $C$ and  $e:\cO_{C}^{n+1}\to\cL$ is
a morphism of sheaves  such
that the Cartier
divisors
$ ( e( \bar{t}_i)) = \sum_{j=1}^{l^i} D_{i,j}. $
\end{enumerate}
With the notations of Definition 1.5 and Proposition 1.6, let $V(\beta)$ be the  subset of $\ol{M}_{g,\cA'}$ representing curves with points $ \{q_{i,j} \}_{0 \leq i \leq n, 1 \leq j
\leq d}$ and $\{p_i\}_{1\leq i\leq m}$ such that $q_{i,k}\not=q_{i,k'}$ whenever $k\in B^i_j$ and $k'\in B^i_{j'}$ with $j\not=j'$.
Let $U({\beta }, \bar{t}):=V(\beta)\times_{\ol{M}_{g,\cA'}}\ol{M}(\bar{t})$ be the $(\CC^*)^n$-- bundle over $V(\beta)$. Then it is not hard to see that $S_{\beta}=S_{b^0}\times ...\times S_{b^n}$ acts as a small group (a group generated by pseudo-reflections) on $U({\beta }, \bar{t})$, and the stabilizer of any point $x\in U({\beta }, \bar{t})$ acts as identity on the fiber $\cC_x$ of the universal family on $\ol{M}(\bar{t})$. Thus to the moduli problem $\ol{\cM}(\mu_{\beta },\bar{t} )$  there corresponds a smooth Deligne-Mumford stack $\ol{M}(\mu_{\beta }, \bar{t}):=U({\beta }, \bar{t})/ S_{\beta}$.


 We show that the map from the Deligne-Mumford stack $U:=\bigsqcup_{\beta, \bar{t}}\ol{M}(\mu_{\beta }, \bar{t})$  to the stack $\ol{M}$ is an  \'etale cover. Indeed, for any $(\cA, a)$-stable family
$( \pi \colon C \to S , \{p_i\}_{1\leq i\leq m}, \cL , e ) $, one can decompose $S=\bigcup_{\beta, \bar{t}}S(\mu_{\beta }, \bar{t})$, where $S(\mu_{\beta }, \bar{t})\subset S$ is the largest open set such that, for any geometric point $s\in S(\mu_{\beta }, \bar{t})$,  the divisor $ ( e( \bar{t}_i))$ on $C_{s}$ can be split into a sum of effective divisors bearing properties (1) and (2). The choice of contributing partitions $\mu_{\beta }$ depends on $a$.

For any $k>0$ we denote by $C^{|k|}$ the $k$-th symmetric product of $C$ over $S$. For any $i=0,...,n$, the Cartier divisor $ ( e( \bar{t}_i))$ mapping to $S$ induces a natural map of $S$ into $C^{|d|}$. The following claim holds.
\begin{claim}
The fiber product $S\times_{\ol{M}}\bigsqcup_{\beta, \bar{t}}\ol{M}(\mu_{\beta }, \bar{t})$ is represented by $\bigsqcup_{\beta, \bar{t}}S'(\mu_{\beta }, \bar{t})$, where each $S'(\mu_{\beta }, \bar{t}):=S(\mu_{\beta }, \bar{t})\times_{(C^{|d|})^{(n+1)}}\prod_{0\leq i\leq n, 1\leq j\leq l^i} C^{|B^i_j|}$. 

Here $\prod$ denotes the fiber product over $S$ and $C^{|B^i_j|}$ stands for the $|B^i_j|$--th symmetric product of $C$.
\end{claim}

 Indeed, pullback of the family $C$ to $S'(\mu_{\beta }, \bar{t})$ satisfies properties (1) and (2). This induces a natural morphism  $S'(\mu_{\beta }, \bar{t})\to \ol{M}(\mu_{\beta}, \bar{t})$. Furthermore, consider any scheme $T$ mapping into both $S$ and $\ol{M}(\mu_{\beta}, \bar{t})$. Let $C_T:=C\times_S T$ and let $C_T^{\beta, \bar{t}}$ be the pullback to $T$ of the universal family over $\ol{M}(\mu_{\beta}, \bar{t})$, and assume that there is an isomorphism over $T$ 
$$\phi_T: C_T \to C_T^{\beta, \bar{t}}$$
compatible with the extra structures $\cL$, $e$ and $\{p_i\}_i$ in a natural way. Then the pullback divisors $\{\phi_T^*(D^i_j)\}_{i,j}$ on $C_T$ induce a canonical map $T\to S'(\mu_{\beta }, \bar{t})$ over $S$. This finishes the proof of the claim. 

 Furthermore, the map $S'(\mu_{\beta }, \bar{t})\to S(\mu_{\beta }, \bar{t})$ is \'etale: $ S(\mu_{\beta }, \bar{t})$ is the quotient of $ S(\mu_{\beta }, \bar{t})\times_{C^{|d|(n+1)}}C^{d(n+1)}$ by the action of the group $(S_d)^{n+1}$, and all  stabilizers are contained in $S_{\beta }$. $S'(\mu_{\beta }, \bar{t})$ is the quotient of $ S(\mu_{\beta }, \bar{t})\times_{C^{|d|(n+1)}}C^{d(n+1)}$ by $S_{\beta }$. This proves that $U$ is an \'etale cover of $\ol{M}$. From here it also follows that the diagonal $\ol{M}\to \ol{M}\times\ol{M}$ is representable.
  
We also notice that it is enough to consider $U$ as  disjoint union of finitely many terms $\ol{M}(\mu_{\beta }, \bar{t})$, and then the fiber product $S'= \bigsqcup_{\beta, \bar{t}}S'(\mu_{\beta }, \bar{t})$ introduced in the Claim  above is of finite type over $S$.
Moreover, given two weighted $(\cA, a)$ -- stable maps
$$\alpha=( \pi \colon C \to S , \{p_i\}_{1\leq i\leq m}, \cL , e ), \alpha'=( \pi' \colon C' \to S , \{p'_i\}_{1\leq i\leq m}, \cL' , e' ),$$
there is a  closed embedding 
$$Isom_S(\alpha, \alpha')\times_SS'(\mu_{\beta }, \bar{t})\hookrightarrow Isom_{S'(\mu_{\beta }, \bar{t})}(\alpha_{S'(\mu_{\beta }, \bar{t})}, \alpha'_{S'(\mu_{\beta }, \bar{t})})\times ((S_d)^{n+1}/S_{\beta }),$$
where $\alpha_{S'(\mu_{\beta }, \bar{t})}, \alpha'_{S'(\mu_{\beta }, \bar{t})}$ are the pullbacks of $\alpha$ and $ \alpha'$ to $S'(\mu_{\beta }, \bar{t})$, together with divisors $ \{D_{i,j}\}_{0\leq i\leq n, 1\leq j\leq l^i}$, and $ \{D'_{i,j}\}_{0\leq i\leq n, 1\leq j\leq l^i}$  as in the definition of $\ol{\cM}(\mu_{\beta }, \bar{t})$. As
$\ol{M}(\mu_{\beta }, \bar{t})$ is a Deligne-Mumford stack, $$Isom_{S'(\mu_{\beta }, \bar{t})}(\alpha_{S'(\mu_{\beta }, \bar{t})}, \alpha'_{S'(\mu_{\beta }, \bar{t})})\to S'(\mu_{\beta }, \bar{t})$$ is separated, and therefore 
$Isom_S(\alpha, \alpha')$ is separated over $S$, which shows that the morphism $\ol{M}\to\ol{M}\times \ol{M}$ is separated. This concludes the proof that  $\ol{M}$ is a Deligne-Mumford stack. 

 Moreover, Hassett in \cite{hassett} proved that the diagonal morphisms for his spaces of wighted curves is finite, which by construction of  $S'(\mu_{\beta }, \bar{t})$ also implies that  $Isom_{S'(\mu_{\beta }, \bar{t})}(\alpha_{S'(\mu_{\beta }, \bar{t})}, \alpha'_{S'(\mu_{\beta }, \bar{t})})$ is finite over $S'(\mu_{\beta }, \bar{t})$. Thus $Isom_S(\alpha, \alpha')$ is finite over $S$, the diagonal morphism $\ol{M}\to\ol{M}\times \ol{M}$ is proper, and $\ol{M}$ is separated. As a particular case of the next Lemma, there is a surjection from $ \ol{M}_{g,m}(\PP^n,d)$ into  $ \ol{M}$ (see also the beginning of section 5), and thus $ \ol{M}$ is proper.

\end{proof}

 We note that if char $k =p$ is positive, and if a weighted stable map factors through a Frobenious morphism $x\to x^p$, then no appropriate $\bar{t}$ cover may be found such that condition (1) in Definition 1.5 is satisfied. Hence the condition that char $k >d$ or char $k =0$ is necessary for Proposition 1.7.


 The next lemma plays an important role in unraveling the structure of the spaces $ \ol{M}_{g,\cA}(\PP^n,d,a)$.

\begin{lemma} Let $m$ and $m'$ be two natural numbers such that $1\leq m'\leq m$. Consider a collection $\cA'$ of $m'$ weights
 and a collection of $m$ weights  $\cA=(\cA', k_1a, ..., k_{m-m'}a)$, where $k_1,..., k_{m-m'}$ are natural numbers such that $\sum_{i=1}^{m-m'}k_i< 1/a$.

There is a natural transformation from the moduli problem of $ (\cA,a) $-stable degree $d$ maps to
$\PP^n$ to the moduli problem of  $
(\cA', a) $-stable,  degree
$d+\sum_{i}k_i$   maps to $\PP^n$. For any $a'<a$, there is a
commutative diagram:

\bea \begin{CD} \ol{M}_{g,\cA}(\PP^n,d,a) @>>>
\ol{M}_{g,\cA'}(\PP^n,d+\sum_{i}k_i, a)\\
@VVV  @VVV \\  \ol{M}_{g,\cA}(\PP^n,d,a') @>>>
\ol{M}_{g,\cA'}(\PP^n,d+\sum_{i}k_i, a')
\end{CD}. \eea

The horizontal maps are local regular imbeddings, while the vertical maps are birational contractions. 

\end{lemma}

\begin{proof}
 Consider an $ (\cA, a) $-stable degree $d$ family of maps to $\PP^n$:
$( \pi :C \to S , \{p_i\}_{1\leq i\leq m}, \cL , e,\mu )$. Endowed
with the line bundle $\cL'=\cL\otimes\cO_C(\sum_i k_i p_i)$ and
the composition
$$\cO_C^{n+1}\to \cL \to \cL\otimes\cO_C(\sum_i k_i p_i),$$
it becomes an  $ (\cA', a)
$-stable degree $d+\sum_{i}k_i$ family of maps to $\PP^n$.

Conversely, given an
$(\cA', a)$-stable degree
$d+\sum_{i}k_i$ family $( \pi :C \to S , \{p_i\}_{1\leq i\leq m'},
\cL' , e,\mu )$ of maps to $\PP^n$ such that the morphism
$\cO^{n+1}\to \cL'$ factors through $\cO^{n+1}\to \cL$ and
provided that $\cL'=\cL\otimes\cO_C(\sum_i k_i s_i)$ for some
sections $(p_i)_{i=m'+1,...,m}$ of $\pi$, we have readily
available  an $ (\cA,a) $-stable degree $d$ family of maps.


The vertical morphisms are induced locally on the rigid covers from Hassett's birational reductive morphisms. These morphisms on the rigid $\bar{t}$-- covers are $(S_d)^{n+1}$-- equivariant. Moreover, the induced maps on the $(S_d)^{n+1}$-- quotients have a natural functorial description that insures their gluing. 
\end{proof}

 Among the spaces of weighted stable curves defined by Hassett, a special role is played by the moduli spaces introduced in \cite{losev}. Their corresponding moduli spaces of weighted stable maps are the simplest spaces on which Gromov-Witten computations can take place, as they are endowed with evaluation maps for all marked points.

\section{The moduli space $\ol{M}_{0,m}(\PP^n,d)$ and its canonical stratification}

 The moduli space  $\ol{M}_{0,m}(\PP^n,d)$ is endowed with a canonical stratification by strata indexed by stable trees, which will play an essential role in our computations. The references here are \cite{behrend1}, \cite{km}.


\begin{definition}


   A graph $\tau$ is given by the data $(F_{\tau}, V_{\tau}, j_{\tau}, \partial_{\tau})$, where $F_{\tau}$ and $V_{\tau}$ are two finite sets with a map $\partial_{\tau} : F_{\tau}\to V_{\tau}$, called the incidence map, and $j_{\tau}: F_{\tau} \to F_{\tau}$ is an involution.

 $F_{\tau}$ is called the set of flags, $V_{\tau}$ the set of vertices , $L_{\tau}=\{f\in  F_{\tau} | j_{\tau}f=f\}$ the set of leaves and $E_{\tau}=\{ \{f_1, f_2\} \subseteq F_{\tau} | f_2 = j_{\tau}f_1 \}$ the set of edges of ${\tau}$. The names are motivated by the geometric realization of  ${\tau}$ as a 1--dimensional topological space, having the vertices as 0--cells, the edges as 1--cells joining two vertices, and the leaves as boundary, each joined with the  vertex corresponding to it via $\partial_{\tau}$ by one of the remaining 1--cells. We will often identify graphs with their geometrical realization.

 For each vertex $v\in V_{\tau}$, the valence $n(v)=|\partial_{\tau}^{-1}(v)|$ counts the half-edges (edges or leaves) adjoint to $v$.
 \end{definition}

 \begin{definition}  A tree $\tau$ is a connected graph such that $|F_{\tau}|=|V_{\tau}|+|E_{\tau}|+|L_{\tau}|-1$.

Set $M=\{1,...,m\}$ and let $D$ be a set of cardinality $d$, where
$m>0$ and $d>0$. An $(M,d)$-tree is a tree $\tau$ along with a
bijection $b: L_{\tau}\to M$ and a function $d: V_{\tau} \to \NN$
such that $\sum_{v\in V_{\tau}} d(v) = d$.

An $(M,D)$-tree is a tree $\tau$  along with a bijection $b:
L_{\tau}\to M\sqcup D$ and a function $d: V_{\tau} \to \NN$ as
above, such that $d(v)$ is exactly the number of leaves in
$L_{\tau}(v)$ labeled by elements of $D$. There is a forgetful map from $(M,D)$-trees to $(M,d)$-trees, which forgets the leaves labeled by $D$. 

 An $(M,d)$-tree $\tau$ is stable if for any vertex $v\in V_{\tau}$, $n(v)>2$ or $d(v)>0$. An $(M,D)$-tree $\tau$ is stable if its associated  $(M,d)$-tree is  stable. The set of $(M,d)$ trees is finite (see \cite{getzlerpan}).

\end{definition}

 We recall here a few basic facts about contractions of trees.

 \begin{definition}
 A contraction of trees $c:\sigma \to \tau$ is given by a pair $(c_V,c_F)$ of a surjective map $c_V: V_{\sigma}\to V_{\tau}$ and an injection $c_F: F_{\tau}\to F_{\sigma}$ which sends $F_{\tau}$ bijectively into $F_{\sigma}\backslash\{f\in F_{\sigma} \backslash L_{\sigma} ;  c_V\partial_{\sigma}(f)=c_V\partial_{\sigma}j_{\sigma}(f) \}$, such that $c_F$, $c_V$ are compatible with $\partial_{\sigma}$, $\partial_{\tau}$,  $j_{\sigma}$,   $j_{\tau}$ and such that, for any vertex $v'\in V_{\tau}$:
  $$F_{\tau}(v')=\bigcup_{v\in c_V^{-1}(v')} F_{\sigma}(v)\backslash\{f\not\in L_{\sigma}(v)  ;  c_V\partial_{\sigma}(f)=c_V\partial_{\sigma}j_{\sigma}(f)=v' \}.$$
Additionally, the geometrical realization of the subset of flags $$\{f\in
F_{\sigma} ; c_V\partial_{\sigma}(f)=v' \}$$ and their incident
vertices is connected.

 If both $c_V$ and $c_F$ are bijections, then $c$ is called an isomorphism. 

 \end{definition}

\begin{definition}
Let $A$ be a semigroup. An $A$-structure on a tree $\tau$ is a
function $\alpha : V_{\tau} \to A$.  A contraction of trees  with $A$-structures $c:(\sigma, \alpha_{\sigma}) \to (\tau,  \alpha_{\tau})$  satisfies the extra condition that for any $v'\in V_{\tau}$:
$$\alpha_{\tau}(v')=\sum_{v\in c_V^{-1}(v')}\alpha_{\sigma}(v).$$

 As such, a contraction of $(M,d)$ or $(M,D)$-trees preserves both the $d$-structure and the labeling of leaves for each vertex.

  \end{definition}

 A contraction also induces injective maps $E_{\tau}\to E_{\sigma}$ and $L_{\tau}\to L_{\sigma}$ in a natural way. If $|E_{\sigma}|=|E_{\tau}|+1$, then $c$ is called a 1-edge contraction. Any contraction can be written in a non-unique way as a composition of 1-edge contractions and automorphisms. Two trees which resulted from the same tree $\tau$ after contraction of the same edges are isomorphic in a unique way.

 We obtain a category $\Gamma_{(M,d)}$ of  $(M,d)$-trees with contractions and final object given by the tree of one vertex of degree $d$ and with $m$ labeled leaves, and a category $\Gamma_{(M,D)}$ of  $(M,D)$-trees with contractions and final object given by the tree of one vertex of degree $d$ and with $m+d$ labeled leaves. There is a natural forgetful functor $\Gamma_{(M,D)} \to \Gamma_{(M,d)}$.

\begin{remark}
We note that the objects in the category $\Gamma_{(M,d)}$ admit automorphisms
other than the identity, whereas the only automorphisms in the
category $\Gamma_{(M,D)}$ are identities. Indeed, assume   $a:\tau
\to \tau$ is an automorphism of $(M,D)$-trees and that $a_V(v)=v'$
for some vertices $v\not=v'$. There is a unique no-return path
connecting $v$ and $v'$, with ending flags $f$ and $f'$ such that
$\partial_{\tau}(f)=v$, $\partial_{\tau}(f')=v'$. By stability,
the branch $\beta_v$ made of all the flags and vertices connected
to $v$ by a  no-return path which does not pass through $f$ must
contain at least one of the labeled leaves $f_0$. Since any
automorphism fixes the labels, $a_T(f_0)=f_0$. Then $a$ must send the unique no-return path from
$f_0$ to $v$ into the no-return path from $f_0$ to $v'$. These
paths differ in their finite number of edges, so there cannot be
an automorphism $a$ interchanging $v$ and $v'$.
\end{remark}


 Each stable map $f: ( C, q_1,..., q_m)\to \PP^n$ of degree $d$ corresponds to a stable $(M,d)$-tree $\tau$, such that 
\begin{itemize}
\item there are bijections between the set of irreducible components of $C$ and the vertices of $\tau$, between the nodes of $C$ and the edges of $\tau$, and between the marked points of $C$ and the leaves of $\tau$,  
\item  these bijections send the incidence relations among components and special points of $C$ (markings or nodes) into corresponding incidence relations among vertices, edges and leaves. 
\end{itemize}
Here $d(v)$ denotes the degree of the component corresponding to $v$.

  Each $(M,d)$-tree $\tau$ uniquely determines a stratum $M(\tau)$ of $\ol{M}_{0,m}(\PP^n,d)$, whose closure $\ol{M}({\tau})$ is the image of the natural morphism $$\phi_{\tau}: \ol{M}_{\tau}\to \ol{M}_{0,m}(\PP^n,d).$$
 Here $\phi_{\tau}$ is an $|Aut(\tau)|$-- degree   ramified cover over $\ol{M}(\tau)$ and $\ol{M}_{\tau}$ denotes the fibered product  ${\prod}^{ev}_{v\in V_{\tau}}  \ol{M}_{0,F_{\tau}(v)}(\PP^n,d(v))$ of moduli spaces along the evaluation maps at the nodes defined by the edges of
$\tau$:
$$\ol{M}_{\tau}\cong {\prod}_{v\in V_{\tau}}  \ol{M}_{0,F_{\tau}(v)}(\PP^n,d(v))\times_{(\PP^n)^{F_{\tau}}} (\PP^n)^{E_{\tau}\sqcup L_{\tau}}$$
 via the natural product of evaluation maps $${\prod}_{v\in V_{\tau}}  \ol{M}_{0,F_{\tau}(v)}(\PP^n,d(v))\to (\PP^n)^{F_{\tau}}=(\PP^n)^{\cup_{v\in V_{\tau}}F_{\tau}(v)}$$ and via the natural map $(\PP^n)^{E_{\tau}\sqcup L_{\tau}}\to (\PP^n)^{F_{\tau}}$  ( see \cite{behrend1}).  $\ol{M}_{\tau}$ parametrizes stable curves along with a contraction of their dual tree to $\tau$. By abuse of terminology, we will sometimes call $\ol{M}_{\tau}$  a normal substratum. 

 There is a natural functor from $\Gamma_{(M,d)}$ to the category of smooth Deligne-Mumford stacks, sending  $\tau$ to $\ol{M}_{\tau}$.
For any contraction  $c:\sigma \to \tau$ of $(M,d)$-trees, the
natural morphism between spaces $\ol{M}_{\sigma}$ and
$\ol{M}_{\tau}$ is a product of gluing morphisms  $
\prod^{ev}_{\{v ; c_V(v)=v'\}}
\ol{M}_{0,F_{\sigma}(v)}(\PP^n,d(v)))  \to
\ol{M}_{0,F_{\tau}(v')}(\PP^n,d(v)))$.
 We are interested in the composition of the functor above with the forgetful functor $\Gamma_{(M,D)}\to\Gamma_{(M,d)}$.


 By analogy with \cite{km} and \cite{losev}, a useful description of $\Gamma_{(M,D)}$ can be given in terms of good monomials. This is based on a notion of a good family of partitions. A detailed analysis of the combinatorics of partitions and trees is done in \cite{losev} for painted stable trees.
The $(M,D)$-trees in this paper closely resemble Losev and Manin's
painted trees. Although the stability condition is different,
the results on good partitions and trees found in \cite{losev} and
\cite{km} also hold for  $(M,D)$-trees.


\begin{definition}

Let $D':=D\sqcup M$.
 A 2-partition $\sigma$ of $D'$ is a pair of two subsets $A_1$, $A_2$ whose union is $D'$.  $\sigma$  is called stable if $A_i\cap D\not=\emptyset$ or $|A_i|\geq 2$ for each $i$. Thus  stable 2-partitions $\sigma =\{A_1, A_2 \}$ of $D'$  are associated to 1-edge stable $(M,D)$-trees, dual to the nodal curves of a boundary divisor in $\ol{M}_{0,m}(\PP^n,d)$. Let $L_{\sigma}$ denote the 1-edge tree determined by the partition $\sigma$.

 For  any stable $(M,D)$-tree $\tau$  there is a family of stable 2-partitions $\{\sigma_e ; e\in E_{\tau}\}$, where $\sigma_e$ is the partition associated   to the tree obtained by contracting all edges except $e$. Thus one can define a monomial $$m(\tau)=\prod_{e\in E_{\tau}} L_{\sigma_e}.$$

 Two stable 2-partitions $\{A_1, A_2\}$ and $\{B_1,B_2\}$ are called compatible if there are exactly 3 non-empty pairwise distinct sets among $A_i\bigcap B_j$, where $i,j= 1$ or 2. A family of stable 2-partitions is called good if any two of them are compatible. The empty family and the family of one stable 2-partition are considered good, too.

 The language of good monomials is particularly well adapted to the description of the Chow ring of $\ol{M}_{0,m}(\PP^n,d)$. Indeed, the divisors $L_{\sigma}$ define an extension  $B^*(\ol{M}_{0,m}(\PP^n,d))$  of the Chow ring, which admits an explicit presentation (conform sections 4 and 5). Classes of the boundary strata  in  $\ol{M}_{0,m}(\PP^n,d)$ are polynomials in $L_{\sigma}$, and only good monomials are non-zero.

\end{definition}

 The following is an analogue of  Lemma 1.2 in \cite{km} and Lemma 1.2.1 in \cite{losev}:

\begin{lemma}
 There is a contravariant equivalence of categories between the category $\Gamma_{(M,D)}$  of stable $(M,D)$-trees and that of good families of 2-partitions of $D'$ whose morphisms are inclusions. At the level of objects, this equivalence is given by   $$\tau \to \{\sigma_e ; e\in E_{\tau}\}.$$
\end{lemma}


\begin{proof}

The bijection at the level of objects is proved in \cite{km}. (Suitable changes in the definition of tree stability do not change the argument). The
tree of $1$ vertex corresponds to the empty family of good
partitions. 

At the level of morphisms, we start with an inclusion
$E\subset E'$ of good families of 2-partitions such that
$E'\backslash E = \{\sigma_1,...,\sigma_l\}$. Let $\tau$ and
$\tau'$ be the corresponding stable $(M,D)$-trees, let
$e_1,...e_l$ be the edges of $\tau'$ associated to
$\sigma_1,...,\sigma_l$ and let $\tau''$ be obtained by
contraction of the edges $e_1,...e_l$. For any other edge $e$ of
$\tau'$, contraction of $\tau'$ to $e$ factors through with the
contractions of $e_1,...e_l$ and thus  $L_{\sigma_e}$ divides both
$m(\tau)$ and $m(\tau'')$. It follows that $m(\tau'')=m(\tau)$ and thus
$\tau= \tau''$ (see Remark 2.5).

\end{proof}

   As a corollary,  there is at most a unique morphism between any two objects of $\Gamma_{(M,D)}$, corresponding to an inclusion of good families of 2-partitions, or equivalently, to divisibility of good monomials. Furthermore, for any two stable $(M,D)$-trees $\tau_A$ and $\tau_B$, there is another stable $(M,D)$-tree $\tau_{A\cap B}$, the final object with respect to contractions of $\tau_A$ and $\tau_B$. There may also exist  a stable $(M,D)$-tree $\tau_{A\cup B}$,  the initial object with respect to contractions to $\tau_A$ and $\tau_B$.
  The corresponding morphisms between fibered products of moduli spaces:
  \bean \diagram
 & {\ol{M}_{A}} \drto^{\phi_A^{A\cap B}} \\
{\ol{M}_{A\cup B}} \urto^{\phi_{A\cup B}^{A}}  \drto_{\phi_{A\cup B}^{B}} & &{\ol{M}_{A\cap B}} \\
 & {\ol{M}_B} \urto_{\phi_B^{A\cap B}} \enddiagram \eean
 and their composition $\phi_{A\cup B}^{A\cap B}= \phi_A^{A\cap B}\circ \phi_{A\cup B}^{A}$ play a central role in the construction of the ring $B^*(\ol{M}_{0,m}(\PP^n,d))$ (section 4). If  $\tau_{A\cup B}$ does not exist, then $\ol{M}_{A\cup B}$ denotes the empty set. 

\section{Closed substrata of the intermediate moduli spaces}

 As with $\ol{M}_{0,m}(\PP^n,d)$, there are categories of trees and of normal strata for each moduli space $\ol{M}_{0,\cA}(\PP^n, d,a)$. For any tree $\tau$ and any vertex $v\in V_{\tau}$, we denote by $E_{\tau}(v)$ the set of edges of $\tau$ incident to $v$.

\begin{definition}

 Let $a\in \QQ$, $a>0$ and let $\cA=(a_1,...,a_m)$ be a collection of weights.
 An $(\cA, a, D, \cR)$-tree is an $(M,D)$ tree $\tau$ along with:
\begin{enumerate}
\item an $\cR$ -- structure, which by definition is a function on $V_{\tau}$ associating to each vertex $v\in V_{\tau}$ a symmetric and transitive relation $\cR(v)\subset L_{\tau}(v)\times  L_{\tau}(v)$, reflexive on the leaves labeled by elements of $M$;
\item  a weight structure $\cA(v)=((a_o)_{o\in O_{\cR(v)}}, (1)_{f\in E_{\tau}(v)})$, where $ O_{\cR(v)}$ is the set of classes (orbits) of $\cR(v)$, and $a_o=\sum_{i\in o} a_i$, with $a_i:= a$ if $i\in D$.
\end{enumerate}
\end{definition}

 Consider an  $(\cA, a)$-- stable map $f:(C, q_1,..., q_m)\to \PP^n$ and its dual $(M,d)$--tree $\tau$ describing the splitting type of $C$. Consider also any  $(M,D)$ tree which is send to  $\tau$ by the natural forgetful functor.
In view of Definition 1.1, the $\cR$-- structure on the tree corresponds to marked points on $C$, or points in the support of $\Coker e$, signaling when these special points coincide.

\begin{definition}
Let $\tau$ be an $(\cA, a, D,\cR)$-tree. Then $a(v):=\sum_{i\in
L_{\tau}(v)} a_i + |E_{\tau}(v)|$.

  The tree  $\tau$ is stable if $a(v)>2$ for any vertex $v$ ($(\cA, a, D)$ stability ) and $\sum_{l\in o} a_l \leq 1$ for any class $o$ of $\cR_{\tau}(v)$ ($\cR$ stability).
\end{definition}

\begin{example}
  When $\cA=(1,...,1)$ and $a>1$, we recover  the usual notion of $(M,D)$-tree stability.
 \end{example}

\begin{definition}
 A contraction  of  $(\cA, a, D, \cR)$-trees preserves the $\cA$, $\cR$ and $d$-structures and the labeling of leaves.

 Definition 2.6 may be adapted  by working with $(\cA,a,D,\cR)$ -- stable 2-partitions of $M\bigsqcup D$. This leads to a notion of good monomials in this context.  As before, contractions between two stable objects $\sigma$ and $\tau$ are unique, can be written uniquely as a division of two good monomials and will contribute to the construction of the ring $B^*(\ol{M}_{0,\cA}(\PP^n,d,a))$. Still, to account for all the natural substrata of $\ol{M}_{0,\cA}(\PP^n,d,a)$ a notion of semi-contraction is also required.

A semi-contraction $c:\sigma\to \tau$  preserves the weight $\cA$, the degree $d$  and
the labeling of leaves  and semi-preserves the equivalence
relation: $\cR_{\tau}(v')\subseteq \bigcup_{v\in c^{-1}_V(v')}
\cR_{\sigma}(v).$ A semi-isomorphism $e:\sigma\to \tau$ is a
semicontraction which is identity at the level of $(M,D)$-trees
with $\cA$-structure, and which semi-preserves $\cR$. Any
semi-contraction may be written in a non-unique way as a
composition of contractions and semi-isomorphisms.
Semi-contractions preserve stability.
\end{definition}


\begin{lemma}
Let $c_{\tau}:\tau'\to \tau $ be  a contraction and $s:
\sigma\to\tau $  a semi-contraction of $(\cA, a, D, \cR)$-trees.
Then there is a unique  $(\cA, a, D, \cR)$-tree $\sigma'$ together
with a contraction $c_{\sigma}$ and a semi-contraction $s'$ making
the following into a pullback diagram:
        \bea \begin{CD}  \sigma' @>{s' }>> \tau'\\
@VV{c_{\sigma} }V  @VV{c_{\tau} }V \\
\sigma @>{s }>> \tau  \end{CD} \eea
  \end{lemma}
\begin{proof}
If $s$ is a contraction, this follows after decomposition into
good monomials. Thus it is enough to check the lemma when $s$ is a
semi-isomorphism. Then $\sigma' = \tau'$ as $(\cA, a, D)$ -- trees and
$c_{\tau }=c_{\sigma }$ as contractions of  $(\cA, a, D)$ -- trees, in the spirit of Definition 2.4.
For $v'\in V_{\sigma'}$ and $v=c_{\sigma V}(v')$,
$\cR_{\sigma'}(v'):=(c_{\sigma F}\times c_{\sigma
F})(\cR_{\sigma}(v))\bigcap (L_{\sigma'}(v')\times
L_{\sigma'}(v'))$. $\cR$-stability along with the universality
property of $\sigma'$ follow accordingly.
\end{proof}

We denote by  $\Gamma_{M,D,\cR}^{\cA, a}$ the category of $(\cA,
a, D,\cR)$-stable trees with semi-contractions.

To each $(\cA,a)$-stable degree  $d$ map to $\PP^n$ there is an
associated stable   $(\cA, a, D, \cR)$-tree. For each  $(\cA, a, D, \cR)$-stable tree $\tau$ there
is a smooth Deligne-Mumford stack $\ol{M}_{\tau}$ parametrizing
$(\cA,a)$-stable degree  $d$ maps along with a semi-contraction of
their dual tree to $\tau$. Specifically,
$$\ol{M}_{\tau}={\prod}^{ev}_{v\in V_{\tau}}
\ol{M}_{0,\cA_{\tau}(v)}(\PP^n,d'(v),a),$$ where ${\prod}^{ev}$
denotes the fibered product along the evaluation maps at the nodes
indexed by $E_{\tau}$, and $d'(v):=d(v)-|\{ f\in L_{\tau}(v)\cap D
| [f]_{\cR}\not=\emptyset\}|$. To a tree contraction there
corresponds a product of gluing maps as usual, and for a
semi-isomorphism $\tau\to\tau'$ sending $v\in V_{\tau}$ to $v'$
there is the map $$\ol{M}_{0,\cA(v)}(\PP^n,d'(v),a)\to
\ol{M}_{0,\cA(v')}(\PP^n,d'(v'), a)$$ described by the horizontal
lines in Lemma 1.8. The images of the spaces $\ol{M}_{\tau}$ in
$\ol{M}_{0,\cA}(\PP^n, d, a)$ are closed strata in a natural
stratification.

We define a notion of stabilization of  $(\cA, a, D, \cR)$-trees
that behaves well with respect to semi-contraction, leading to a
functor $\Gamma_{M,D,\cR}^{\cA, a} \to \Gamma_{M,D,\cR}^{\cA',
a'}$ for appropriate choices of $(\cA, a)$ and $(\cA', a')$, as described in Theorem 1.4.

\begin{definition}

Let $\tau$ be an $(\cA, a, D, \cR)$-tree, $v_0\in V_{\tau}$ and
$f\in F_{\tau}(v_0)$.  The branch $\beta(f)$ made of all the flags
and vertices connected to $v_0$ by a  no-return path  passing
through $f$ is called unstable if the sum $\sum_{v\in
V_{\beta(f)}\backslash\{v_0\}} a(v) \leq 2| E_{\beta(f)}|.$ 
 If $\tau$ is dual to a curve $C$, and $C'\subset C$ corresponds to the branch $\beta(f)$, then by the inequality above $C$ is $(\cA,a)$-- unstable, because condition (1) of Definition 1.1 fails on $C'$. 

A stable morphism $z: \tau\to \tau'$ is a contraction of $(\cA, a,
D)$ trees which contracts  all edges of unstable branches
$\beta(f)$, and such that \bea &  \cR_{\tau'}(z_V(v_0))= (z_F^{-1}\times
z_F^{-1})\cR_{\tau}(v_0)\bigcup &\\ & \bigcup((\cup_{v\in
V_{\beta(f)}\backslash\{v_0\}}z^{-1}_FL_{\tau}(v))\times (\cup_{v\in
V_{\beta(f)}\backslash\{v_0\}}z^{-1}_FL_{\tau}(v))).&\eea For any
$(\cA, a, D, \cR)$-tree $\tau$ there exists a unique $(\cA, a, D,
\cR)$-stable tree $\tau_s$ and stable morphism $s: \tau\to \tau_s
$ into it, such that $s$ factors through any other stable morphism of
source $\tau$. The stabilization $s$ contracts all the edges of
all the unstable branches of $\tau$ and $\tau_s$ is $(\cA, a, D,
\cR)$-stable.

\end{definition}

\begin{lemma}
 For any semi-contraction of $(\cA, a, D, \cR)$-trees $c: \sigma
\to \tau$, there is a unique semi-contraction of $(\cA, a, D,
\cR)$- stable trees $c_s:\sigma_s \to \tau_s $ making the
following into a pullback diagram:
        \bea \begin{CD}  \sigma @>{c }>> \tau\\
@VV{s_{\sigma} }V  @VV{s_{\tau} }V \\
\sigma_s @>{c_s }>> \tau_s  \end{CD}. \eea
\end{lemma}

\begin{proof}
The lemma is immediate if $c$ is a semi-isomorphism. In
that case $c_s$ is also a semi-isomorphism. Thus it is enough to
check the above in the case when $c$ is a 1-edge contraction. Let
$e$ be the contracted edge. There are two cases: when the edge is
part of an unstable branch and when it is not. In the second case,
$c_s$ is simply the contraction of the edge $s_{\tau}(e)$. In the
former case $c_s$ is clearly a semi-isomorphism.

\end{proof}

\begin{remark}

In particular, given two collections of $m$ weights $\cA$ and
$\cA'$ such that $a_i\geq a'_i$, and given $a>a'$, to the diagram of Lemma 3.7 there
corresponds a cartesian diagram
   \bea \begin{CD}  \ol{M}_{\sigma} @>>> \ol{M}_{\tau}\\
@VV{s_{\sigma}}V  @VV{s_{\tau} }V \\
\ol{M}_{\sigma'} @>>> \ol{M}_{\tau'} \end{CD}, \eea
 where the objects on the upper horizontal line map to closed substrata
of $\ol{M}_{0,\cA}(\PP^n,d,a)$ and those on the lower line map to
closed substrata of $\ol{M}_{0,\cA'}(\PP^n,d,a')$.
\end{remark}

\section{The extended Chow ring $B^*(\ol{M}_{0,\cA}(\PP^n,d,a))$}

There is a natural action of the group of symmetries $G:=S_d$ on the
category of $(\cA, a, D, \cR)$--trees, induced by permutatations in the set of leaves marked by elements of $D$. Similarly, $G$ acts on the set of 2-partitions as well as on the set of good families of  2-partitions  associated to $(\cA, a, D, \cR)$--trees. For any good family $P$ of 2-partitions, let $G_P\subset G$ denote its stabilizer, which fixes all elements of $P$. 

 The moduli spaces $\ol{M}_P$ corresponding to good families of $(\cA,a, D, \cR)$-stable 2-partitions of $M\sqcup D$ form a network of regular local embeddings, with commutative diagrams like the one at the end of section 2. For every $[ g ]\in G/G_P$ there is a canonical isomorphism $g: \ol{M}_P\to \ol{M}_{g (P)}$, induced by the isomorphism of $(\cA, a, D, \cR)$--trees. Moreover, from the $(\mu_{\beta },\bar{t})$--covers described in the proof of Proposition 1.7 we  may extract   a set of \'etale covers $\ol{M}_P(t)\to \ol{M}_P$, such that
there is a Cartesian diagram
\bea  \diagram  \bigsqcup_{[ g ]\in G_N/G_P} \ol{M}_{g(P)}(t)
\rto \dto & \ol{M}_N(t)\dto \\
\ol{M}_P \rto^{{\phi}_P^N } & \ol{M}_N,
\enddiagram \eea
each morphism $\ol{M}_P(t)\to \ol{M}_N(t)$ is an inclusion and $\ol{M}_P(t)\bigcap \ol{M}_Q(t)=\ol{M}_{P\bigcup Q}(t)$ in $\ol{M}_{P\bigcap Q}(t)$, and all intersections are transverse. 

 The local \'etale structure of this network, and the action of the group $G$ on it, allow the construction of an extended Chow ring for  $\ol{M}_{0,\cA}(\PP^n,d,a)$ by concatenating the Chow rings of all network spaces, and identifying cycles with their pushforwards in a way compatible with the diagram above.   

\begin{definition}

  Given a collection of weights $\cA$ and a positive number $a$, we construct  a graded $A^*(\ol{M}_{0,\cA}(\PP^n,d,a))$- algebra  \bea B^*(\ol{M}_{0,\cA}(\PP^n,d,a)):=     ( \bigoplus_{l=0}^{\dim (\ol{M}_{0,\cA}(\PP^n,d,a))} (\oplus_P A^{l-|P|}(\ol{M}_P)) ) / \sim ,   \eea
where the second sum is taken after all  good families of $(\cA
,a, D, \cR)$-stable 2-partitions of $M\sqcup D$, corresponding to trees
that can be contracted to the tree of 1 vertex, $(d+m)$ -- leaves. The
equivalence relation $\sim$ is generated by
$$\phi^{N}_{P*}(\alpha)\sim \sum_{[g]\in G_N/G_P} g_*(\alpha)$$   
for any two  good families of 2-partitions $N$ and $P$ such that $N\subset P$ and any class $\alpha\in A^*(\ol{M}_P)$.



Due to Lemma 3.13 in \cite{noi}, this definition coincides to Definition 3.6 in \cite{noi}.

  For any two  good families  of 2-partitions $A$ and $B$, and any two classes $\alpha \in A^*(\ol{M}_A)$ and $\beta \in
A^*(\ol{M}_B)$,  multiplication is defined by:
$$ \alpha\cdot \beta := \phi^{A*}_{A\cup B}(\alpha) \cdot \phi^{B*}_{A\cup
B}(\beta) \cdot c_{top}(\phi^{A\cap B*}_{A\cup
B}\cN_{\ol{M}_{A\cap B} | \ol{M}_{\emptyset}})$$ in
$A^*(\ol{M}_{A\cup B})$. Here we keep notations from diagram (2.1). $\phi^{A*}_{A\cup B},
\phi^{B*}_{A\cup B}$ are the generalized Gysin homomorphisms
for regular local embeddings, as defined in \cite{vistoli}.
$\cN_{\ol{M}_{A\cap B} | \ol{M}_{\emptyset}}$ is the normal bundle
(as stacks) of $\ol{M}_{A\cap B}$ in $\ol{M}_{\emptyset}:=
\ol{M}_{0,\cA}(\PP^n,d,a)$.

\end{definition}

 The multiplication is well defined. The proof is identical to that of Lemma 3.12 in \cite{noi}.

Given any $(\cA, a, D, \cR)$- stable tree $\tau$, the 
$ A^*(\ol{M}_{\tau})$-- algebra $B^*(\ol{M}_{\tau})$ may be defined analogously to Definition 4.1, by considering only trees that admit contractions to $\tau$. Let $P$ be the good family of 2-- partitions associated to $\tau$. The group $G_P$ admits a natural action on $B^*(\ol{M}_{\tau})$. There is naturally defined Reynolds operator $\Psi: B^*(\ol{M}_{\tau})\to A^*(\ol{M}_{\tau})$ and by the definition of $\sim$ we obtain the following 

\begin{lemma}

The subring of invariants of $B^*(\ol{M}_{\tau})$
under the action of $G_P$ is $A^*(\ol{M}_{\tau} )$. 
\end{lemma}

In particular, the subring of invariants of $B^*(\ol{M}_{0,\cA}(\PP^n,d,a))$
under the action of $S_d$ is $A^*(\ol{M}_{0,\cA}(\PP^n,d,a))$.

For any semi-contraction of trees $\sigma\to \tau$ and the associated gluing morphism of moduli spaces $\phi : \ol{M}_{\sigma}\to \ol{M}_{\tau}$, there is a natural pullback  $\phi^*: B^*(\ol{M}_{\tau})\to B^*(\ol{M}_{\sigma})$, defined via the diagram of Lemma 3.5.
Indeed, from definition it follows that $\sim$ is compatible with the morphisms
induced by semi-contractions as well as stabilizations of trees.


Given two collections of $m$ weights $\cA$ and $\cA'$ such that
$a_i\geq a'_i$, and $a>a'$, to the morphism
$f:\ol{M}_{0,\cA}(\PP^n,d,a)\to \ol{M}_{0,\cA'}(\PP^n,d,a')$ there
is associated a pullback morphism
$f^*:B^*(\ol{M}_{0,\cA'}(\PP^n,d,a')) \to
B^*(\ol{M}_{0,\cA}(\PP^n,d,a))$ constructed in the following way.
For any $(\cA, a, D, \cR)$- stable tree $\sigma$ admitting a contraction
to the tree $\iota $ of 1 vertex and $d+m$ leaves, there is a
unique $(\cA', a', D, \cR)$- stable tree $\sigma'$ and a cartesian diagram
 \bea \begin{CD}  \ol{M}_{\sigma} @>>> \ol{M}_{0,\cA}(\PP^n,d,a)\\
@VV{s_{\sigma}}V  @VV{s }V \\
\ol{M}_{\sigma'} @>>> \ol{M}_{0,\cA'}(\PP^n,d,a') \end{CD}, \eea
induced by the diagram of Lemma 3.7.

Let  $\iota' $ be the $(\cA', a', D, \cR)$-- stable tree of 1 vertex and $d+m$ leaves.
The semi-contraction $\sigma' \to \iota'$ can be written uniquely
as a composition of a semi-isomorphism $i:\sigma'\to \sigma''$
with a contraction. The  semi-isomorphism $i$ differs from identity precisely when the map $s$ is a weighted blow-down and $\ol{M}_{\sigma}$ is (the normalization of) an exceptional divisor in $\ol{M}_{\sigma''}\times_{\ol{M}_{0,\cA'}(\PP^n,d,a')}\ol{M}_{0,\cA}(\PP^n,d,a)$.  The morphism $f^*$ is
obtained by summing up the pullbacks $
A^*(\ol{M}_{\sigma''})\to A^*(\ol{M}_{\sigma})$ for all $\sigma$ as above.

 More generally, there exists a pullback morphism $f_{\tau}^*: B^*(\ol{M}_{\tau_s})\to  B^*(\ol{M}_{\tau})$ for any $(\cA, a, D, \cR)$-- stable tree $\tau$ and its $(\cA', a', D, \cR)$-- stabilization $\tau_s$, constructed as above.




\section{The structure of the ring $B^*(\ol{M}_{0,m}(\PP^n,d))$}

 From here throughout $m>0$.
 
 For integers $k=0,...,d-1$, let $\ol{M}_{0,m}(\PP^n,d,a_k)$ be the moduli space of $(\cA_k, a_k)$-stable maps, where $a_k=\frac{1}{k+\epsilon }$, and the collection of weights of the $m $ marked points is $\cA_k:=(1, a_k, ..., a_k)$ if $k>0$, and $\cA_0:=(1,1,...,1)$. In particular, $\ol{M}_{0,m}(\PP^n,d,a_0)=\ol{M}_{0,m}(\PP^n,d)$.


As a special case of  Lemma 1.8, for $k\geq 1$ there is a commutative diagram

\bea \begin{CD} \ol{M}_{0,m}(\PP^n,d,a_k) @>>>  \ol{M}_{0,1}(\PP^n,d+m-1,a_k)\\ @VVV  @VVV \\    \ol{M}_{0,m}(\PP^n,d,a_{k+1})        @>>> \ol{M}_{0,1}(\PP^n,d+m-1,a_{k+1})  \end{CD}. \eea

The morphisms $ \ol{M}_{0,1}(\PP^n,d+m-1,a_k) \to \ol{M}_{0,1}(\PP^n,d+m-1,a_{k+1})$ have been described in detail in \cite{noi} as weighted blow-downs. The imbedding 
$$\ol{M}_{0,m}(\PP^n,d,a_k)\to\ol{M}_{0,1}(\PP^n,d+m-1,a_k)$$ establishes $\ol{M}_{0,m}(\PP^n,d,a_k)$ as one of the natural normal strata of $\ol{M}_{0,1}(\PP^n,d+m-1,a_k)$. This is self-evident at the level of trees. Indeed, let $M=\{1_M,2_M,...,m_M\}, D=\{1_D,...d_D\}$, $M'=M\backslash \{1_M\}$, $D'=M'\sqcup D$ and $d'=|D'|=d+m-1$. Any tree $\tau\in \Gamma^{\cA_k,a_k}_{M,D,\cR}$ may be thought of as an element of $ \Gamma^{(1),a_k}_{\{ 1_M\}, D',\cR'}$ by simply switching the leaves pertaining to $M'$ from the role of marked-point holders to the role of degree holders. $\tau$ may be regarded as rooted tree by taking as root $r$ the vertex to which the leaf $1_M$ is attached. 

In particular, $\ol{M}_{0,m}(\PP^n,d,a_k)=\ol{M}_{\iota_m}$ where $\iota_m\in  \Gamma^{(1),a_k}_{\{ 1_M\}, D',\cR'}$ is the tree with $1$ vertex and $d+m$ leaves, such that $\cR(\tau )= \{(i_M, i_M); i=1,...,m \}.$ In \cite{noi}, section 3, $ \ol{M}_{0,1}(\PP^n,d',a_k)$ is denoted by $ \ol{M}_{0,1}(\PP^n, d', d'-k)=\ol{M}^{d'-k}$, while  $\ol{M}_{0,m}(\PP^n,d,a_k)$ is $\ol{M}^{d'-k}_{\{\{2_M\},...,\{m_M\}\} }$ for $k\geq 1$.

 $\ol{M}_{0,m}(\PP^n,d)$ is obtained by successive weighted blow-ups along strata given by trees semi-isomorphic to $\iota_m$. The effects of these transformations on Chow rings are described in \cite{noi}, Lemmas 3.15 to 3.21. In the present context we  obtain an analogue to \cite{noi}, Proposition 3.22:

 The algebra $B^*(\ol{M}_{0,m}(\PP^n,d))$ is generated over $\QQ$ by codimension 1 classes $$\{T_{\sigma }\}_{\sigma }, H, \mbox{ and } \psi,$$  where ${\sigma }$ ranges over all $(M,D)$-stable 2-partitions of the set $M\sqcup D$. $H$ is the pull-back of the hyperplane class in $\PP^n$ by the evaluation map at the point $1_M$, and $\psi$ is the tautological class at the same point. Let $\sigma =(h_1,h_2)$ be represented by the set $h_1$ such that $1_M\not\in h_1$. We will write $T_{h_1}$ for $T_{\sigma }$, as in \cite{noi}.  $-T_{h_1}$ represents the fundamental class of the divisor $$\ol{M}_{\sigma}=\ol{M}_{h_1\bigcap M, \star}(\PP^n, |h_1\bigcap D|)\times_{\PP^n}\ol{M}_{\star, h_2\bigcap M}(\PP^n, |h_2\bigcap D|), $$ which is pullback of the divisor class for $\ol{M}_{0,1}(\PP^n, |h_1|)\times_{\PP^n}\ol{M}_{0,2}(\PP^n, |h_2|-1)$  from $B^*(\ol{M}_{0,1}(\PP^n,d'))$. We will denote this class also by $-T_{h_1}$. 



 For any stable tree $\tau\in \Gamma^{(1),a_k}_{\{ 1_M\}, D',\cR'}$ with $1\leq k\leq d'-1$ and any edge $e\in E_{\tau}$, there is a relation in $B^*(\ol{M}_{\tau})$ whose geometric significance is explained below.

Think of $\tau$ as a rooted tree. We say that a vertex $w'$ is subordinated to another vertex $w$ if any path from $w'$ to $r$ passes through $w$. For any $w\in V_{\tau}$, let $\beta_{w}$ denote the branch of $w$, which is the rooted subtree of root $w$ containing all vertices, edges and leaves subordinated to $w$. The degree $d_w$ of $\beta_w$ is the cardinal of the set of leaves $ L_{\beta_w}$. With the notations above, $T_w:=T_{L_{\beta_w}}$. Thus in  $B^*(\ol{M}_{0,1}(\PP^n,d'))$, $-T_w$ is the fundamental class of $\ol{M}_{w}:=\ol{M}_{0,1}(\PP^n,d_w)\times_{\PP^n}\ol{M}_{0,2}(\PP^n,d'-d_w)$.  Set
 \bean H_w:=H+(d'-d_w)\psi +\sum_{L_{\beta_w}\subset h\subset D'} (|h|-d_w)
 T_h,\eean
 \bean \psi_w=\psi+ \sum_{L_{\beta_w}\subseteq h\subset D'}
 T_h.\eean
 $H_w$ is the pullback of the hyperplane class by the evaluation map, and $\psi_w$ the canonical class, at the marked point of  $\ol{M}_{0,1}(\PP^n,d_w)$, pulled back to $\ol{M}_{w}$. (see \cite{ leepan}).

We define the set $$A_w:=   \{w'\in V_{\beta_w}  ; w' \mbox{ immediately subordinated to } w\}. $$ Let $e=(u,v)\in E_{\tau}$ be an edge of vertices $u$ and $v$, with $v$ subordinated to $u$. Consider the partial partitions of $D'$
\bean I=\{ L_{\beta_w} | {w\in A_v} \}\bigcup \cR'[v],  J=\{ L_{\beta_w} | w\in A_u\backslash \{v\} \} \bigcup \cR'[u],\eean
 where for any vertex $w$, $\cR'[w]$ is the set of all classes of elements in $L_w$ given by  the relation $\cR'$. In particular, when $\tau$ comes from an  $(M,D)$-tree, then  $\cR'[w]\supseteq M\cap L_w.$

 Let $\tau_1$ denote the tree obtained from $\tau$ after contraction of the edge $e$ and $a_{d_v}$-stabilization, and let $\tau_2$ be the tree obtained by $a_{d_v}$-stabilizing $\tau$, such  that $\tau_2$ is quasi-isomorphic to $\tau_1$. (Note that fixing weight $a_{d_v}$ on the leaves labeled by $D$ makes the entire branch $\beta(v)$ unstable.) The following relation satisfied by the ``exceptional divisors'' at the weighted blow-up of $\ol{M}_{\tau_1}$ along the image of $\ol{M}_{\tau_2}$, is the pullback of relation (4), Proposition 3.22 in \cite{noi} applied to $\ol{M}_{0,1}(\PP^n, d_u)$:
\bean  \left. \prod_{w\in V_{\tau_1}} T_{w}\left(\sum_{h\in
H(I,J)} P_{\tau,e}(\{T_{h'}\}_{h'\supset h},
t_{h})\right|^{t_{h}=T_{h}}_{t_{h}=0}+ P_{\tau,e}(0)\right) \eean
We will denote this relation by $r(\tau,e)$.
 Here $H(I,J)=\{h\subseteq L_{\beta_u}  ;  |h|\geq d_v, h\supseteq h_I \mbox{ and } h\cap h_J=\emptyset \mbox{ whenever } |h|=d_v \}$, where $h_I:=\bigcup_{h'\in I}h'$, $h_J:=\bigcup_{h'\in J}h'$   and for $\emptyset\not=h\in H(I,J),$
$P_{\tau,e}$ is the equivalent of $P^{d'-k}_{I,J}$ defined in \cite{noi}, (3.7).\bea P_{\tau,e}(\{t_{h'}\}_{h'\supseteq h})= \psi'^{i-1}
(\prod_{j=1+|h|-d_v}^{|h\backslash h_I|}(H'+j\psi')^{n+1}-
H'^{(n+1)(d_v-|h_I|)}), \eea where $i:= |A_v|+|\cR'[v]|$ and
$$H'=H_u+(d_u-|h\cup h_J|)\psi_u+\sum_{ L_{\beta_u}\supset h'\supset h} |h'\backslash (h\cup h_J )|T_{h'},$$
 $$\psi'=\psi_u+\sum_{ L_{\beta_u}\supset h'\supseteq h} T_{h'}. $$
When $h=\emptyset$, we denote by $$P_{\tau,e}(0):=\psi_u^{i-1}
\prod_{j=1+d_u-d_v-|J|}^{d_u-|I\cup J|}(H_u+j\psi_u)^{n+1}. $$

 When $k=d'-1$, the stable trees $\tau\in \Gamma^{(1),a_{d'}}_{\{ 1_M\}, D',\cR'}$ have no edges, and correspond to the simplest moduli spaces and substrata, before any weighted blow-ups have been performed. By comparison with \cite{noi}, Proposition 3.22, (3) \cite{noi}, pullback of relations from the extended Chow ring of these spaces may be encoded as $r(\tau, e)$, where $e$ is the pair made of two copies of the root. 

 Other natural relations in $B^*(\ol{M}_{0,m}(\PP^n,d);\QQ)$ are given by the compatibility conditions: $\prod_{\sigma\in A} T_{\sigma} =0$ whenever $A$ is not a set of good monomials. Finally, there is the relation $H^{n+1}=0$, pulled back from $\PP^n$ by the evaluation morphism at the first marked point. The equivalent relations for the other marked points are included among (5.4). These, together with $r(\tau, e)$ for all trees $\tau$ coming from stable $(M,D)$-trees, and all edges $e$, form a complete set of relations for the ring $B^*(\ol{M}_{0,m}(\PP^n,d);\QQ)$, due to Lemmas 3.15 and 3.19 in \cite{noi}, and by an  identical argument to that employed in the proof of Proposition 3.22, \cite{noi}.
 However, the Appendix to \cite{noi} shows that the above set of relations is redundant. Indeed, any relation introduced at the $a_k$-th intermediate step of weighted blow-down is in the ideal generated by relations coming up at the $a_{k-1}$-th intermediate step. A complete and less redundant set of relations can be formulated as follows:

\begin{theorem}
 Let $M=\{1_M, 2_M,...,m_M\}, D=\{1_D,...d_D\}$, $M'=M\backslash \{1_M\}$, $D'=M'\sqcup D$ and $d'=|D'|=d+m-1$.

$B^*(\ol{M}_{0,m}(\PP^n,d);\QQ)$ is the $\QQ$-algebra generated by divisors $$H \mbox{ , } \psi \mbox{ , } T_{h}$$ for all $ h\subset D'$ such that $h\not=\emptyset$ or $\{ i_M\}$ for $i_M\in M'$. Let $T_{\emptyset}:=1$.

The ideal of relations is generated by:

\begin{enumerate}
\item $H^{n+1}$;
\item $T_{h}T_{h'}$ unless $h\cap h'=\emptyset$, or $\emptyset\not= h\subseteq h'$ or $\emptyset\not= h'\subseteq h$;
\item \begin{itemize}

\item ($m\geq 1$) $T_{h}T_{h'}(\psi+\sum_{h\cup h'\subseteq h''} T_{h''})$ for all $h \not= h'$ nonempty;
\item ($m \geq 2$) $T_{h}(\psi+\sum_{h\cup\{i_M\}\subseteq h'} T_{h'})$ for all  $h\not=\emptyset$ and  $i_M\in M'\backslash h $;
\item ($m\geq 3$) $\psi + \sum_{\{i_M, j_M\}\subseteq h} T_{h}$ for all $i_M, j_M \in M' $;
\end{itemize}

\item $(m>1)$ $(H+d\psi +\sum_{i_M\in h} |h\cap M'| T_{h}))^{n+1}$ for all  $i_M\in M'$;


\item $T_{h}(\sum_{h'\not=h}\left. P(t_{h'})\right|^{t_{h'}=T_{h'}}_{t_{h'}=0}+ \psi^{-1}(H+|^ch_D|\psi)^{n+1})$ for all $h$,

where \bea  && P(t_{h'})= (\psi+ \sum_{h'' \supset h'} T_{h''} +t_{h'})^{-1} \\
&& [ (H+|^ch_D| \psi +\sum_{h'' \supset h'}|h_D''\backslash h_D| T_{h''} +|h_D'\backslash h_D| t_{h'})^{n+1}-\\
&&-(H+|^ch_D\cap ^ch_D')|\psi +\sum_{ h'' \supset h'  }|h_D''\backslash (h_D\cup h_D')| T_{h''})^{n+1}].  \eea

Here for any $h \subset D'$, $h_D:=h\cap D$ and $^ch_D:=D\setminus h$.
\end{enumerate}

\end{theorem}

\begin{proof}

 Keeping the same notations as before, our task is to eliminate the redundant relations among $r(\tau,e)$.
  We notice that $r(\tau,e)$ is
  actually pullback of relation $r(\sigma,e )$, where $\sigma$ is
  the stable $(M,D)$-tree obtained after contracting all edges of
  $\tau$ not adjacent to $u$ or $v$, such that $V_{\sigma}=\{r,u,v\}\cup A_u\cup
  A_v$. Furthermore, by Remark 4.2 in \cite{noi},
  $$r(\sigma,e)=T_u( r(\sigma',e)-r(\sigma'',e)),$$
  where $\sigma'$ is obtained after contracting the edge $(r,u)$
  of $\sigma$ and $\sigma''$ is obtained from $\sigma$ after replacing the flag
  $(e,u)$ by a flag $(e,r)$. Thus we may consider relations
  corresponding to edges adjacent to the root of the tree only. Now Lemma
  4.3 in \cite{noi} implies that all relations existing in an
  $a_k$--th intermediate space are combinations of relations appearing, after a weighted
  blow-up, at the $a_{k-1}$--th intermediate space. Recall the sets $I$, $J$
  introduced above. Here $I$, $J$ are associated to an edge $(r,v)$ and 
 $\cR'[v], \cR'[r] \subseteq M'$.
Following Lemma 4.4 in \cite{noi}, it is enough to consider the
cases when $|I|\leq 2, |J|\leq 1$ such that $|I|$ and $J$ are
minimal as to insure the stability of $\sigma$. Thus the minimal
building blocks for all relations in $B^*(\ol{M}_{0,m}(\PP^n,d))$
besides (1) and (2) are relations based on trees satisfying one of the following
\begin{itemize}
\item $d(v)=0$ and $I =L_{{w_1}}\cup
L_{{w_2}}$ or $I =L_{{w}} \cup \{i_M\}$ or
 $I =\{i_M, j_M \}$. Clearly here one can restrict to the case when $J =\emptyset$, as  the stability requirement is satisfied.
 The ensuing relations are (3).
\item $d(v)=1$, $I =\{i_M\}$,  $J =\emptyset$.
 This generates relation (4).
Note that this relation
cannot be further reduced because $T_{i_M}$ does not exist. If,
however, one considers $I =L_{{w}}$, the
ensuing relation is reduced by Lemma 4.4 in \cite{noi} to the
following:
\item $d(v)=1$, $I =\emptyset$, $J =\emptyset$ or
$J =L_{{w}}$. These are relations (5).

\end{itemize}

\end{proof}

Relation (3) $(m\geq 3)$ is  the classically known
formula for the $\psi$-class on $\ol{M}_{0,m}(\PP^n,d)$.
 Relation (4) is pullback of the relation on $\PP^n$
by the evaluation map at $i_M$. Due to the choice of the element $1_M\in M$, the structure of the ring as presented above is not symmetrical in all marked points. In particular the well-known divisorial relations obtained by pullback from $\ol{M}_{0,4}$ do not appear here because some of their terms are already not among the ring generators listed by us. 

 The  relations analogous to (1)-(5) obtained when fixing a different privileged marked point in $M$ can be checked to be dependent of (1)-(5) via Remarks 4.1 and 4.2 in \cite{noi} and due to the divisorial relations mentioned above. 

\section{An additive basis of $A^*(\ol{M}_{0,m}(\PP^n,d))$} 

We can index the elements of an additive basis of
$A^k(\ol{M}_{0,m}(\PP^n,d))$ by automorphism classes of rooted
$(M,d)$-trees $\tau$ with a weight function on vertices $b:
V_{\tau}\to \NN$, satisfying  the following conditions:
\begin{itemize}
\item The root $r$ is the vertex to which the leaf $1_M$ is attached. It satisfies $$0 \leq b(r) < (n+1)d(r)+ n(r)-2 .$$
\item For any vertex $v\not= r$, $$ 0 < b(v) < (n+1)d(v)+ n(v)-2 .$$
\item The sum of weights $$ \sum_{v\in V_{\tau}} b(v) \leq k.$$
\end{itemize}
 A tree with the above extra-structure will be called an $(M,d,b,k)$-tree. If  the last condition is omitted, the resulting structure is called an $(M,d,b)$-tree.
 As usual, an automorphism of a rooted $(M,d,b)$-tree is an automorphism of the underlying tree, which fixes the labeling of its leaves and root, and preserves the functions $d(v)$ and $b(v)$.

 To an automorphism class $[\tau]$ of rooted $(M,d,b,k)$-trees  we associate the class $[\tau_k] \in A^k(\ol{M}_{0,m}(\PP^n,d))$ as follows: First, we choose any partition $p: V_{\tau} \to \cP(D)$ compatible with the degree map $d(v)$. By abuse of notation, let  $\tau$ also denote the $(M,D)$-tree obtained by attaching to the vertices of the $(M,d)$-tree $\tau$ leaves labeled according to $p$. Lemma 2.7 shows how to associate to the $(M,D)$-tree $\tau$ a unique good monomial
   $$M(\tau):=\prod_{e\in E_{\tau}} T_{\sigma(e)},$$
where $\sigma(e)$ is the unique partition of the set $M\bigsqcup D$ obtained after contracting all edges of $\tau$ except $e$. Since $\tau$ is rooted, for any edge $e$ there exists a unique flag $f=(v,e)$ such that the root $r$ can be reached from $v$ by a no-return path starting at $f$. Let $T_v:=T_{\sigma(e)}$. Define \bea [\tau_k] = H^{k-\sum_{v\in V_{\tau }} b(v)}\psi^{b(r)}\mbox{ sym}_d(\prod_{v\in V_{\tau}\backslash\{ r\}} T_{v}^{b(v)}),\eea
where sym$_d$ is the symmetrizing function with respect to the action of the group $S_d$ on the set of good monomials, induced by its action on $D$: $\mbox{ sym}_d (M)= \frac{1}{d !} \sum_{\mu\in S_d} \mu (M)$.

\begin{proposition}
The classes $[\tau_k]$ corresponding to all $(M,d,b,k)$-trees  form
an additive basis for the group $ A^k(\ol{M}_{0,m}(\PP^n,d))$.
\end{proposition}

\begin{proof}
 We have seen that $H, \psi $ and $T_A$ generate $ B^*(\ol{M}_{0,m}(\PP^n,d))$
 as a $\QQ$-algebra, when $A$ ranges over all $(M,D)$-stable two-partitions
 of the set
  $D'=(M\backslash\{1_M\})\sqcup D$. Thus monomials
 $ H^{k-\sum_{v\in V_{\tau }} b(v)}\psi^{b(r)}\mbox{ sym}_d(\prod_{v\in V_{\tau }\backslash\{ r\}} T_{v}^{b(v)})$, for all $(M,D)$-trees $\tau$ and all weight functions $b(v)$ with $k\geq \sum_{v\in V_{\tau }} b(v)$, generate $ A^k(\ol{M}_{0,m}(\PP^n,d))$ as a vector space over $\QQ$.
Furthermore, from equation (5.4) it follows that the classes
$[\tau_k]$ associated to rooted $(M,d,b,k)$-trees are generators for $
A^k(\ol{M}_{0,m}(\PP^n,d))$, since the degree of $P_{\tau,e}$ is
exactly $(n+1)d(v)+n(v)-2$, the codimension of the blow-up locus
$\ol{M}_{\tau_2}$ in $\ol{M}_{\tau_1}$.
 The linear independence is proved by induction on the
 intermediate step $a_l$, decreasing in $l$, and increasing
 induction on dimension of normal strata. The first step is provided by Lemma 3.3 in \cite{noi}.
 Consider now a relation at step $a_l$ on a closed normal
 stratum $\ol{M}'$. Such a normal stratum appears in Lemma 3.19 in \cite{noi} under the notation $\ol{M}'=\ol{M}_I^{k+1}$, where $k=d-l-1$. According to that lemma, $B^*(\ol{M}')$  is generated over $\QQ$ by $ H, \psi,$ and $\{T_{h'}\}_{|h'|>l}$. Thus for any relation among the generators of $B^*(\ol{M}')$, one may choose a set $h$ of minimal cardinality $|h|> l$, and write the relation as $$\sum_{j=0}^{b} c_j
 T_{h}^j=0,$$  where each coefficient $c_j\in
 \QQ[H, \psi, \{T_{h'}\}_{|h'|\geq |h|, h\not= h'\subset D'}]$. 
Let $\ol{M}_h$ be
 the normal stratum of  $\ol{M}_{0, m}(\PP^n, d,
 a_{|h|})$ corresponding to the partition $(h, D\setminus h)$.
 The relation above must also exist in the image $\ol{M}$ of $\ol
{M}'$ in $\ol{M}_{0, m}(\PP^n, d, a_{|h|-1})$. Assume that $b$ is
less than the codimension
 of $\ol{M}\times_{\ol{M}_{0, m}(\PP^n, d, a_{|h|})}\ol{M}_h$ in
 $\ol{M}$. Then by Lemma 3.16 in
 \cite{noi}, each $c_j$ must be in the kernel of the pullback
 morphism $B^*(\ol{M})\to B^*(\ol{M}\times_{\ol{M}_{0, m}(\PP^n, d,
 a_{|h|})}\ol{M}_h)$. Thus the linear independence problem is
 reduced to a smaller substratum and possibly step index $l$.

\end{proof}

  Proposition 6.1 leads to an alternate way to \cite{getzlerpan} of computing the Betti numbers of $\ol{M}_{0,m}(\PP^n,d)$, by tracking down the contribution to the Poincair\'e polynomial of each stable $(M,d)$-tree, given by all the possible $b$- structures on it. However in lack of a closed formula, this computation becomes cumbersome for large $d$. It is more efficient when the number $m$ is concerned.

 First we fix some notations: Let $$P_{\ol{M}_{0,m}(\PP^n,d)}(q):=\sum_{k}\dim_{\QQ } A^k(\ol{M}_{0,m}(\PP^n,d)) q^k.$$

The following corollary follows directly from Proposition 6.1:

\begin{corollary}

The following formula holds:
$$P_{\ol{M}_{0,m}(\PP^n,d)}(q)= \frac{q^{n+1}-1}{q-1}\sum_{(\tau,b)/iso}\prod_{v\in V_{\tau}}q^{b(v)}$$
where the sum is taken after all isomorphism classes $(\tau,b)/iso$  of $(M,d,b)$-trees. 

\end{corollary}

The factor $\frac{q^{n+1}-1}{q-1}$ is the contribution of the pullback of the 
hyperplane section of $\PP^n$ via the evaluation morphism
$ev_1$.

For a rooted stable $(M,d)$-tree $\tau$ without isomorphisms, the contribution to $P_{\ol{M}_{0,m}(\PP^n,d)}(q)$ brought by all possible $b$-structures on it is $$P_{\tau}=q^{|V_{\tau}|-1}\prod_{v\in V_{\tau}} \frac{q^{ (n+1)d(v)+ n'(v)-2}-1}{q-1},$$
 where $n'(v)$ is the number of directly subordinated flags: $n'(r)=n(r)$ and $n'(v)=n(v)-1$ if $v\not=r$. The more involved part of the computation comes from the trees with automorphisms.

\begin{example}

  In the following, let $(k)$ denote a vertex of degree $k$ and let $i_M$ be a leaf corresponding to the $i$--th marked point. With our conventions, the following trees contribute to $P_{\ol{M}_{0,2}(\PP^n,2)}$:


\begin{enumerate}

\item  $\xy \xymatrix{ & {(2)}  \POS[];[dl]**\dir{-}, [dr]**\dir{-}\\
                        {1_M} & & {2_M}  } \endxy$   \hskip  0.2in   $\frac{q^{2(n+1)}-1}{q-1}$;

\item  $\xy \xymatrix{ & {(0)} \POS[]; [dl]**\dir{-}, [];[dr]**\dir{-},  [];[d]**\dir{-}\\
                     {(2)} &   {1_M} & {2_M}  } \endxy$   \hskip  0.1in    $\frac{q(q^{2n}-1)}{q-1};$

\item   $\xy \xymatrix{ & {(1)} \POS[]; [dl]**\dir{-}, [];[dr]**\dir{-},  [];[d]**\dir{-}\\
                     {(1)} &   {1_M} & {2_M}  } \endxy$  \hskip  0.1in    $\frac{(q^{n+2}-1)q(q^{n-1}-1)}{(q-1)^2};$

\item   $\xy \xymatrix{ & {(0)} \POS[]; [dl]**\dir{-}, [];[dr]**\dir{-},  [];[d]**\dir{-}\\
                     {(1)} \POS[]; [d]**\dir{-}, &   {1_M} & {2_M} \\  {(1)}  } \endxy$      \hskip  0.1in $\frac{q^2(q^{n}-1)(q^{n-1}-1)}{(q-1)^2};$

\item    $\xy \xymatrix{ & {(1)} \POS[]; [dl]**\dir{-}, [];[dr]**\dir{-}\\
                        {(1)} \POS[]; [d]**\dir{-}, & & {1_M}\\  {2_M} } \endxy$     \hskip  0.2in   $\frac{(q^{n+1}-1)q(q^{n}-1)}{(q-1)^2};$

\item   $\xy \xymatrix{ & {(0)} \POS[]; [dl]**\dir{-}, [];[dr]**\dir{-},  [];[d]**\dir{-}\\
                     {(1)} \POS[]; [d]**\dir{-}, &   {(1)} & {1_M} \\  {2_M}  } \endxy$  \hskip  0.1in  $\frac{q^2(q^{n}-1)(q^{n-1}-1)}{(q-1)^2};$

\item      $\xy \xymatrix{ & & {(0)} \POS[];   [];[dll]**\dir{-}, [];[dl]**\dir{-}, [];[d]**\dir{-},  [];[dr]**\dir{-}\\
               {(1)} &  {(1)} &   {1_M} & {2_M}  } \endxy$       \hskip  0.1in  $(q+1)\frac{q^2}{2}\left[\frac{(q^{n-1}-1)^2}{(q-1)^2}+\frac{q^{2(n-1)}-1}{q^2-1}\right];$

\end{enumerate}

 Adding all summands and multiplying by the factor contributed by $\PP^n$, we obtain $P_{\ol{M}_{0,2}(\PP^n,2)}$ as in \cite{cox}:

$$ P_{\ol{M}_{0,2}(\PP^n,2)}(q)=\frac{(q^{n+1}-1)(q^{n}-1)(q^{n+3}-1+2q^{n+2}-2q+2q^{n+1}-2q^2)}{(q-1)^3}.$$

\end{example}

\begin{remark} The formula of Corollary 6.2 is based on attaching the leaf $1_M$ to the root of each tree $\tau$. Similarly we may define a larger class of polynomials: for each $l\leq m$, let
$$P^l_m(d):=\sum_{(\tau',b)/iso}\prod_{v\in V_{\tau'}}q^{b(v)},$$
where the sum is taken after all isomorphism classes of pairs $(\tau', b)$ such that the leaves $1_M,..., l_M$ are attached to the root of $\tau'$. Let $S^lP^0_0(d)$ be defined by the same formula, where the sum is now taken after all  isomorphism classes of pairs $(\tau'', b)$  such that $\tau''$ is a collection of $l$ connected components, all rooted $(M,d)$-trees with no leaves. Finally, for any such polynomials $P$, $Q$, let $P\star Q (d)=\sum_{e=0}^d P(e) Q(d-e)$. With these notations, the following formulas hold for any $d$:
\bean P^{l-1}_m(d) = P^{l}_m(d) +q\sum_{j=0}^{m-l} \left(\begin{array}{c} m-l\\j\end{array}\right) P^0_{j+1}\star P^l_{m-j} (d).\eean
Let $\mu\vdash d$ denote a partition of $d$, namely a set of natural numbers $0<\mu_1\leq...\leq\mu_l$, whose sum is $d$. Let $I_1,...,I_r\subseteq \{1,...,l\}$ be a partition of the set of indexes such that $\mu_i=\mu_j$ for any $i,j\in I_s$ and $\mu_i\not=\mu_j$ for any $i\in I_s$, $j\in I_s'$ such that $s\not=s'$. Let $\mu_{I_s}:=\mu_i$, $i\in I_s$.
\bean P^m_m(d)=P^{n+1+m}_{n+1+m}(d-1)+\sum_{\mu\vdash d}q^l \prod_{s=1}^r S^{|I_s|}P^0_0(\mu_{I_s})P^{m+l-2}_{m+l-2}(0).\eean
 To obtain the first formula, consider trees $\tau'$ which contribute to $P^{l-1}_m(d)$ but not to $ P^{l}_m(d)$ and split each of them into two rooted trees, by erasing the edge which is adjoint to the root and in the path from the root to $l_M$. To obtain the second formula, whenever the root of an $(M,d)$ tree $\tau'$ has degree at least 1, we decrease the degree by 1, while adjoining $n+1$ leaves to it, and whenever the root has degree 0, we split the tree by erasing all edges adjoint to the root. Formulas (6.1) and (6.2) give an inductive algorithm for computing $P^l_m(d)$, by increasing  $d$ and decreasing  $l$, starting from the building blocks  $S^iP^0_0(k)$, ($ki\leq d$). Then
$$P_{\ol{M}_{0,m}(\PP^n,d)} = \frac{q^{n+1}-1}{q-1}P^1_m(d).$$
\end{remark}

\begin{example}

$\bullet$ $S^iP^0_0(1)=\frac{(q^{n-1}-1)(q^n-1)...(q^{n+i-2}-1)}{(q-1)(q^2-1)...(q^{i}-1)}$.

$\bullet$ $P^m_m(1)=\frac{q^m-1}{q-1}\frac{q^n-1}{q-1}$ if $m\geq 1$, and  $\frac{q^{n-1}-1}{q-1}$ if $m=0$.

$\bullet$  $P^m_m(2)=\frac{q^m-1}{q-1}\frac{(q^n-1)(q^{n+1}-1)}{(q-1)(q^2-1)}(q^2+q+1)$, if $m\geq 1$, $\frac{(q^n-1)(q^{n+1})}{(q-1)^2}$ if $m=0$. 

$\bullet$   $P^m_m(3)=\frac{q^m-1}{q-1}\frac{(q^n-1)(q^{n+1}-1)}{(q-1)(q^2-1)}[(q^2+q+1)\frac{q^{n+2}-1}{q-1}+q^2(q+1)\frac{q^{n-1}-1}{q-1}+q^3(q-1)\frac{q^{n-1}-1}{q^3-1}]$ if $m\geq 1$.

\end{example}



\begin{remark}
The system of tautological rings $R^*(\ol{M}_{0,m}(\PP^n,d))$ is
defined as the set of smallest $\QQ$-subalgebras of
$A^*(\ol{M}_{0,m}(\PP^n,d))$, such that
\begin{itemize}
\item The system is closed under pushforwards via the forgetful
maps $\ol{M}_{0,m}(\PP^n,d) \to\ol{M}_{0,m-1}(\PP^n,d)$ and the
natural gluing maps $\ol{M}_{\tau}\to \ol{M}_{0,m}(\PP^n,d)$, for
all stable $(M,d)$-trees $\tau$;
\item The system contains all pullbacks of classes from $\PP^n$ by the evaluation map
 $ev :\ol{M}_{\tau}\to \PP^{n}$ corresponding to any leaf or edge
 of $\tau$.
\end{itemize}
It is known that
$A^*(\ol{M}_{0,m}(\PP^n,d))= R^*(\ol{M}_{0,m}(\PP^n,d))$. See \cite{oprea1}, \cite{oprea2} for two different proofs. With our
notations from Proposition 6.1, we can write the classes $[\tau_k]$
as tautological classes as follows: First, all the canonical
classes $\psi_{i_M}$ are tautological by Lemma 2.2.2 in
\cite{pandharipande2}. Thinking of the stable $(M,D)$-tree $\tau$
as rooted like before, for any vertex $u\in V_{\tau}\setminus
\{r\}$ we let $\psi_u$ denote the cotangent class on
$\ol{M}_{0,F_{\tau}(u)}(\PP^n, d(u))$ at the marked point indexed
by the flag in $F_{\tau}(u)$ which is closest to the root. Then
the pushforward $i_{\tau *}(\psi_u) $ by $i_{\tau}: \ol{M}_{\tau}
\to\ol{M}_{0,m}(\PP^n,d)$ is written in
$B^*(\ol{M}_{0,m}(\PP^n,d))$ as:
$$i_{\tau *}(\psi_u)=\sym_{d}[\psi_u\prod_{v\in  V_{\tau}\setminus
\{r\}} T_v]=\sum_{\mu\in S_d}[\psi_{\mu(u)}\prod_{v\in
V_{\tau}\setminus \{r\}} T_{\mu(v)} ] .$$ On the other hand, by
formula (5.2), \bea    \prod_{v\in  V_{\tau}\setminus \{r\}}
T_v^{b(v)}= \prod_{v\in  V_{\tau}\setminus \{r\}} T_v
(\psi_v-\sum_{h\supset L_{\beta_v}}T_h)^{b(v)-1}. \eea After expanding and applying $\sym_d$,  we may reiterate the procedure as needed, until  no terms $T_h^b$ with $b>1$ appear in the expansion. Thus  $\sym_d(  \prod_{v\in  V_{\tau}\setminus \{r\}}
T_v^{b(v)})$, and consequently $[\tau_k]$ are written as a tautological classes.

\end{remark}


\providecommand{\bysame}{\leavevmode\hbox to3em{\hrulefill}\thinspace}


\begin{thebibliography}{DL}


\bibitem[BaM]{bayer-manin} A. Bayer, Yu. Manin, Stability Conditions, Wall-crossing and weighted Gromov-Witten Invariants in math.AG/0607580

\bibitem[BeM]{behrend1}
K.Behrend, Yu.Manin, Stacks of stable maps and Gromov-Witten
invariants in Duke Math. J.  \textbf{85} (1996), 1-60.

\bibitem[BO]{behrend2} K.Behrend, A.O'Halloran, On the cohomology of
stable map spaces in \emph{Invent. Math.} \textbf{154} (2003), no.
2, 385--450.

\bibitem[C]{cox}
J.Cox, An additive basis for the Chow ring of
$\ol{M}_{0,n}(\PP^r,2)$, in math.AG/0501322

\bibitem[C2]{cox2}
J.Cox, A presentation for the Chow ring $A^*(\bar{M}_{0,2}(P^1,2))$ in
math.AG/0505112

\bibitem[FP]{fultonpan}
W.Fulton, R.Pandharipande, Notes on stable maps and quantum
cohomology in \emph{Algebraic geometry, Santa Cruz 1995}, page
45--96 volume \textbf{62} of Proc. Symp. Pure Math. Amer. Math.
Soc. 1997.

\bibitem[GP]{getzlerpan}
E.Getzler, R.Pandharipande, The Betti numbers of
$\ol{M}_{0,n}(r,d)$, in math.AG/0502525


\bibitem[G]{givental}
A.Givental, Equivariant Gromov-Witten invariants in Internat. Math. Res. Notices (1996), no. 13, 613-663.


\bibitem[Has]{hassett}
B. Hassett, Moduli spaces of weighted pointed stable curves. Adv.
Math. \textbf{173} (2003), no. 2, 316-352.


\bibitem[Ke]{keel}
S.Keel, Intersection theory on moduli space of stable n-pointed
genus zero, trans. A.M.S. 1992 volume \textbf{330} page 545--574

\bibitem[Kn]{knudsen}
F.F.Knudsen, The projectivity of the moduli space of stable
curves, II in Math.Scand. 1983, volume \textbf{52}, page 161--199

\bibitem[KM1]{km}
M.Kontsevich, Yu.Manin, Quantum cohomology of a product. With an
appendix by R. Kaufmann. in Invent. Math. \textbf{124} (1996), no.
   1-3, 313--339.


\bibitem[KM2]{kontsevich}
M.Kontsevich, Yu.Manin, Gromow-Witten classes, quantum cohomology,
and ennumerative geometry, in Comm.Math.Physics \textbf{164}
(1994), page 525-562

\bibitem[LM]{losev}
A.Losev, Yu.Manin,  Extended modular operad  in Frobenius
manifolds,  181--211, Aspects Math., E36, Vieweg, Wiesbaden, 2004.


\bibitem[LLY]{yau}
B.H.Lian, K.Liu, S-T.Yau,  Mirror principle. I. Asian J. Math. 1
(1997), no. 4,page 729--763.

\bibitem[LP]{leepan}
Y.-P.Lee, R.Pandharipande, A reconstruction theorem in quantum
cohomology and quantum K-theory, math.AG/0104084

\bibitem[MM1]{noi}
A.Mustata, A.Mustata, Intermediate moduli spaces of stable maps,
in math.AG/0409569, to appear in Invent. Math.

\bibitem[MM2]{noi3}
A.Mustata, A.Mustata, Universal relations on stable map spaces in genus zero,
 preprint

\bibitem[O1]{oprea1}
D.Oprea, The tautological rings of the moduli spaces of stable
maps, in math.AG/0404280

\bibitem[O2]{oprea2}
D.Oprea, Tautological classes on the moduli spaces of stable maps
to projective spaces, in math.AG/0404284

\bibitem[P]{pandharipande2}
R.Pandharipande, Intersection of $\QQ$-divisors on Kontsevich's
Moduli Space $\ol{M}_{0,n}(\PP^r, d)$ and enumerative geometry,
Trans.Amer.Math.Soc. \textbf{351} (1999), no.4, 1481-1505


\bibitem[V]{vistoli}
A.Vistoli, Intersection theory on algebraic stacks and their
moduli spaces in Inv. Math. 1989, volume \textbf{97}, page
613--670



\end{thebibliography}
\end{document}